\documentclass[12pt]{amsart}
\usepackage[utf8x]{inputenc}
\usepackage{amsmath, amsthm, amssymb,cite,enumerate,slashed}
\usepackage{fullpage}
\usepackage{color}
\usepackage{hyperref}
\usepackage[english]{babel}
\usepackage{framed}

\usepackage{pgfplots}
\usepgfplotslibrary{fillbetween}

\newtheorem{thm}{Theorem}
\newtheorem{lem}{Lemma}[section]
\newtheorem{exm}[lem]{Example}
\newtheorem{prop}[lem]{Proposition}
\newtheorem{cor}[lem]{Corollary}

\newtheorem{rmkk}[lem]{Remark}
\newtheorem{rmk}{Remark}[lem]
\newtheorem*{rmkno}{Remark}
\newtheorem{que}[lem]{Question}

\newtheorem{defn}[lem]{Definition}
\numberwithin{equation}{section}

\newcommand{\R}{\mathbb{R}}
\newcommand{\D}{\mathcal{D}}

\newcommand{\eps}{\varepsilon}
\newcommand{\norm}[1]{\left\lVert#1\right\rVert}
\newcommand{\normw}[1]{\left\lVert#1\right\rVert_{\rm w}}
\newcommand{\mcl}[1]{\mathcal{ #1 }}
\newcommand{\fm}[1]{\[\begin{aligned} #1 \end{aligned}\]}
\newcommand{\eq}[1]{\begin{equation}\begin{aligned} #1 \end{aligned}\end{equation}}
\newcommand{\lra}[1]{\langle #1 \rangle}

\newcommand{\khw}[1]{\left( #1 \right)_{\rm w}}
\newcommand{\wt}[1]{\widetilde{ #1 }}

\DeclareMathOperator{\supp}{supp}

\begin{document}

\title{A uniqueness theorem for 3D semilinear wave equations satisfying the null condition}

\author{Dongxiao Yu}
\address{Mathematical Institute, Hausdorff Center for Mathematics, University of Bonn}
\email{yudx@math.uni-bonn.de}

\begin{abstract}

In this paper, we prove a uniqueness theorem for a system of semilinear wave equations satisfying the null condition in $\mathbb{R}^{1+3}$. Suppose that  two global solutions with $C_c^\infty$ initial data have equal initial data outside a ball and equal radiation fields outside a light cone. We show that these two solutions are equal either outside a hyperboloid or everywhere in the spacetime, depending on the sizes of the ball and the light cone. 

\end{abstract}

\maketitle
\tableofcontents\addtocontents{toc}{\protect\setcounter{tocdepth}{1}}

\section{Introduction}

This paper is devoted to a uniqueness result for a system of semilinear wave equations for $u=(u^I)_{I=1}^N$  in $\R^{1+3}_{t,x}$, of the form
\eq{\label{swe}\Box u^I&=(-\partial_t^2+\Delta_x)u^I=Q^I(u,\partial u),\qquad I=1,2,\dots,N.}
Here $Q=(Q^I)$ is a smooth (i.e.\ $C^\infty$) $\R^N$-valued  function defined in  $\R^{N}\times\R^{4N}$. In addition, we assume that  $Q$ satisfies the null condition. That is, if we have the Taylor expansion \eq{\label{swea1} Q^I(u,v)&=\sum_{J,K=1}^N\sum_{\alpha,\beta=0}^3A^{\alpha\beta}_{I,JK}v_\alpha^J v_\beta^K+O((|u|+|v|)^3),\qquad\text{for all } (u,v)\in\R^N\times\R^{4N},}
then for each $I,J,K$, we have
\eq{\label{swea2}A_{I,JK}^{\alpha\beta}\widehat{\omega}_\alpha \widehat{\omega}_\beta=0,\qquad \text{for all }\widehat{\omega}=(-1,\omega)\in\R\times\mathbb{S}^2.}
Note that in \eqref{swea1}, all the $A_{*,**}^{**}$'s are  constants independent of $u$ and $v$. It is known that this system \eqref{swe}, along with initial data
\eq{\label{init} (u,\partial_tu)|_{t=0}=(\eps u_0,\eps u_1)\in C_c^\infty(\R^3;\R^N)\times C_c^\infty(\R^3;\R^N),} has a global solution for sufficiently small $\eps>0$. This global existence result was proved by Klainerman \cite{klai,klai2} and Christodoulou \cite{chri}. In addition, because of the null condition,  this global solution $u$ behaves as a linear solution to the linear wave equation $\Box w=0$ as time goes to infinity. In addition, each  global solution  has a Friedlander radiation field $F_0(q,\omega)$ with $q=r-t$ and $\omega=x/r$. 

In this paper, we give a uniqueness theorem related to the radiation fields. Let $u,v$ be two global solutions to \eqref{swe} for all $t\geq 0$ with  $C_c^\infty$ initial data at $t=0$. We show that if the initial data of two solutions are equal outside the  ball $\{|x|\leq R_1\}\subset\R^3$, and if the radiation fields of these two solutions are equal in the region $q>R_2$ where $R_2<R_1$, then $u=v$ in a region determined by $R_1$ and $R_2$. If $R_1\leq |R_2|$, this region is the whole spacetime $\R^{1+3}_+=\R^{1+3}\cap\{t\geq 0\}$. If $R_1>|R_2|$, this region is  \fm{t\geq 0,\quad r^2-(t+R_2)^2\geq R_1^2-R_2^2.} In the linear case, i.e.\ if $Q^I\equiv 0$ for each $I$, this region is \fm{t\geq 0,\quad r+t\geq R_1,\quad r-t\geq R_2.}
See Figure \ref{fig} for an illustration of these sets. 

\subsection{Nonlinear wave equations}
Let us consider  a generalization of the system \eqref{swe} in~$\R^{1+3}_{t,x}$
\begin{equation}\label{nlw}
\square u^I=F^I(u,\partial u,\partial^2u),\qquad I=1,2,\dots,N.
\end{equation}
The nonlinear term is assumed to be smooth with the Taylor expansion
\begin{equation}\label{nonlinearity}F^I(u,\partial u,\partial^2u)=\sum a_{\alpha\beta,JK}^I\partial^\alpha u^J\partial^\beta u^K+O(|u|^3+|\partial u|^3+|\partial^2 u|^3).\end{equation}
The sum  is taken over all $1\leq J,K\leq N$ and all multiindices $\alpha,\beta$ with $|\alpha|\leq |\beta|\leq 2$, $|\beta|\geq1$ and $|\alpha|+|\beta|\leq 3$. Besides, the coefficients $a_{\alpha\beta,JK}^I$'s are all universal constants. 

\subsubsection{Lifespan}
Since 1980's, several results  on the lifespan of the solutions to the Cauchy problem \eqref{nlw} with initial data \eqref{init} have been proved. Let us first focus on the scalar case when $N=1$. For example, John \cite{1john,1john2} proved that \eqref{nlw} does not necessarily have a global solution; in fact, any nontrivial solution to  $\square u=u_t\Delta u$ or $\square u=u_t^2$ must blow up in finite time. In contrast, in $\R^{1+d}$ with $d\geq 4$, H\"{o}rmander \cite{1horm3} proved  the small data global existence for \eqref{nlw}. For arbitrary nonlinearities in three space dimensions, the best result on the lifespan is the almost global existence: the solution exists for $t\leq \exp(c/\eps)$ where $\eps\ll 1$. The almost global existence for \eqref{nlw} was proved by Lindblad \cite{1lind3}, and  we also refer to   \cite{1johnklai,1klai,1horm2,1horm,MR2015331,MR1945285} for some earlier work. In the case when $N>1$, if the nonlinear terms $(F^I)$ in \eqref{nlw} do not depend on $u$ (i.e.\ if $F^I=F^I(\partial u,\partial^2 u)$), then we still have small data global existence if $d\geq 4$ and almost global existence if $d=3$. However, if the nonlinear terms depend on $u$, then we would encounter some tricky situations. We refer our readers to \cite{metcmorg,metcrhoa} for two papers which reveal the differences between the scalar case $N=1$ and the general case.

In contrast to the finite-time blowup in John's examples,  Klainerman \cite{klai2} and Christodoulou \cite{chri} proved that the null condition is sufficient for  small data global existence. The null condition, first introduced by Klainerman \cite{klai}, states that for each $1\leq I,J,K\leq N$ and for each $0\leq m\leq n\leq 2$ with $m+n\leq 3$, we have 
\begin{equation}\label{nullcond}A_{mn,JK}^I(\omega):=\sum_{|\alpha|=m,|\beta|=n}a_{\alpha\beta,JK}^I\widehat{\omega}^{\alpha}\widehat{\omega}^{\beta}=0,\hspace{1cm}\text{for all }\widehat{\omega}=(-1,\omega)\in\R\times\mathbb{S}^2.\end{equation}
Equivalently, we assume that $A_{mn,JK}^I\equiv 0$  on the null cone $\{m^{\alpha\beta}\xi_\alpha\xi_\beta=0\}$ where $(m^{\alpha\beta})=(m_{\alpha\beta})$ is the Minkowski metric $\text{diag}(-1,1,1,1)$. The null condition leads to cancellations in the nonlinear terms \eqref{nonlinearity} so that the nonlinear effects of the equations are much weaker than the linear effects. 
We also remark that the null condition is not necessary for  small data global existence. One such example is the Einstein vacuum equations in wave coordinates; see \cite{1lindrodn,1lindrodn2}.  We also refer our readers to \cite{1tata} for a general introduction to the null condition.

Later, Lindblad and Rodnianski \cite{1lindrodn,1lindrodn2} introduced the weak null condition. To state this condition, we start with the asymptotic equations first introduced by H\"{o}rmander \cite{1horm2,1horm,1horm3}. We make the ansatz\begin{equation}\label{ansatz}
u^I(t,x)\approx \eps r^{-1}U^I(s,q,\omega),\hspace{1cm}r=|x|,\ \omega_i=x_i/r,\ s=\eps\ln(t),\ q=r-t,\ 1\leq I\leq N.
\end{equation}
Assuming that $t=r\to\infty$, we substitute this ansatz into \eqref{nlw} and compare the coefficients of  terms of order $\eps^2 t^{-2}$.  Nonrigorously, we can obtain the following asymptotic PDE
\begin{equation}\label{asypde11}2\partial_s\partial_q U^I=\sum A_{mn,JK}^I(\omega)\partial_q^mU^J\partial_q^nU^K.\end{equation}
Here $A_{mn,JK}^I$ is defined in \eqref{nullcond} and the sum is taken over $1\leq J,K\leq N$ and $0\leq m\leq n\leq 2$ with $m+n\leq 3$. We say that the \emph{weak null condition} is satisfied if \eqref{asypde11} has a global solution for all $s\geq 0$ and if the solution and all its derivatives grow at most exponentially in $s$, provided that the initial data decay sufficiently fast in $q$. In the same papers, Lindblad and Rodnianski  conjectured that the weak null condition is sufficient for small data global existence. To the best of the author's knowledge,  this conjecture  still remains open, but we refer to Keir \cite{1keir,1keir2} for some recent progress. For more discussions on the weak null condition, we refer our readers to Section 1.1 of the author's recent paper \cite{yu2021asymptotic}.

\subsubsection{Radiation field}
The radiation field for the wave equation was first studied by Friedlander \cite{fried1,fried2,fried3,fried4}. For simplicity, we start with an $\R$-valued solution $u$ to a linear wave equation $\Box u=0$ with initial data $(u,u_t)|_{t=0}\in C_c^\infty(\R^3)$. Then, it can be proved that the limit
\eq{\label{f0defn}F_0(q,\omega)=\lim_{r\to\infty} ru(r-q,r\omega)}
exists for each $q\in\R$ and $\omega\in\mathbb{S}^2$, and that the function $\R\times\mathbb{S}^2\ni (q,\omega)\mapsto F_0(q,\omega)$ is a smooth function. In addition, we have a pointwise estimate
\eq{\label{f0pointest}|Z^I (u(t,x)-r^{-1}F_0(r-t,\omega))|\lesssim_{I} t^{-2},\qquad \forall I,\ r>t/2>1.}
Here $Z^I$ is a product of $|I|$ vector fields where each $Z$ represents one of the commuting vector fields: translations $\partial_\alpha$, scaling $t\partial_t+r\partial_r$, rotations $x_i\partial_j-x_j\partial_i$ and Lorentz boosts $x_i\partial_t+t\partial_i$; see Section \ref{sec2.1}. Such a function $F_0$ is called the \emph{radiation field}. For the proofs of these results, we refer our readers to Friedlander \cite{fried1,fried2,fried3,fried4} and Section 6.2 of H\"{o}rmander \cite{1horm}.

It turns out that most of the discussion above still applies to  global solutions to \eqref{swe} with small and localized initial data. That is, if $u=(u^I)$ is a global solution to \eqref{swe} for $t\geq 0$ with initial data \eqref{init} for sufficiently small $\eps$, then the limit $F_0(q,\omega)=(F_0^I(q,\omega))$ defined by \eqref{f0defn} exists and is smooth with respect to $(q,\omega)\in \R\times\mathbb{S}^2$. Moreover, for each integer $M>0$, the pointwise estimates \eqref{f0pointest}, with an additional factor $\lra{r-t}$ on the right hand side, hold for all $|I|\leq M$, as long as $\eps\ll_{M}1$. Intuitively, these results remain true because the null condition guarantees that  the nonlinear effects of the system \eqref{swe} are much weaker than the linear effects. For the sake of completeness, we will present brief proofs of these results in Section \ref{sec2.2}.

Inspired by these  examples, we give the following definition which will be used in this paper.

\defn{\label{defrad}\rm Suppose that $u=(u^I)$ is a global solution to \eqref{nlw} for $t\geq 0$ with $C_c^\infty$ initial data. We say that an  $\R^N$-valued $C^1$ function $F_0=(F_0^I)$ of $(q,\omega)\in\R\times\mathbb{S}^2$ is the \emph{radiation field} of $u$, if 
\eq{\label{def1a1}F_0(q,\omega)=\lim_{r\to\infty}ru(r-q,r\omega),\qquad\forall  (q,\omega)\in\R\times\mathbb{S}^2,} and if
\eq{\label{def1a2}\sum_{|I|\leq 1}|Z^I(u(t,x)- r^{-1}F_0(r-t,\omega))|\lesssim t^{-2}\lra{r-t},\qquad\forall r>t/2\gtrsim 1.}
Here $Z^I$ is a product of $|I|$ vector fields where each $Z$ represents one of the commuting vector fields: translations $\partial_\alpha$, scaling $t\partial_t+r\partial_r$, rotations $x_i\partial_j-x_j\partial_i$ and Lorentz boosts $x_i\partial_t+t\partial_i$.
}
\rmk{\label{mrmk1.10}\rm The Lorentz boosts $\Omega_{0i}$ for $i=1,2,3$  in Definition \ref{defrad}  are in fact not necessary in our proof. Note that the commuting vector fields are introduced because we need the pointwise estimate \eqref{l2.1c1} below. However, in this paper, we are only interested in the region where $|r-t|\lesssim 1$ and $t\gtrsim 1$, so instead of  \eqref{l2.1c1} we can apply
\fm{|(\partial_t+\partial_r)\phi|+|(\partial_i-\omega_i\partial_r)\phi|\lesssim r^{-1}(\lra{r-t}|\partial\phi|+|S\phi|+\sum_{1\leq j<k\leq 3}|\Omega_{jk}\phi|).}
We include the Lorenzt boosts in this paper only for simplicity.}
\rmk{\rm A direct corollary of \eqref{def1a2} is that for each fixed $q^0\in\R$, we have \eq{\label{defradcor}\sum_{|I|\leq 1}|Z^Iu|\lesssim_{q^0} (1+t+r)^{-1},\qquad \forall r-q^0>t>0.}This pointwise estimate will be useful in the rest of this paper. To prove \eqref{defradcor}, we first note that $u\equiv 0$ whenever $r-t\geq C$ and that $F_0\equiv0$ whenever $q\geq C$ for some constant $C$.  These two identities follow from the finite speed of propagation. Besides, we note that $|u|+|Zu|\lesssim 1$ for $t\lesssim1$. So from now on, we can assume $q^0\leq r-t\leq C$ and $t\gtrsim 1$. Since $Z^I(r-t)=O(\lra{r-t})$ and $Z^I\omega=O(1)$ whenever $r\sim t$, by the chain rule we have
\fm{\sum_{|I|\leq 1}|Z^I(r^{-1}F_0(r-t,\omega))|\lesssim_{q^0} (1+t+r)^{-1},\qquad \forall |r-t|\lesssim q^0+1,\ t\gtrsim 1.}
This finishes the proof of \eqref{defradcor}.}
\rmk{\rm The radiation fields are introduced to study the asymptotic behavior of the global solutions to \eqref{nlw}. However, it is not guaranteed that a global solution to \eqref{nlw}  admits a radiation field in general. In those cases, we need to introduce some other notion of ``radiation fields''. For example, in \eqref{def1a1}, instead of taking the limit along a straight line, we may take the limit along a characteristic, i.e.\ a null curve with respect to some Lorentzian metric related to \eqref{nlw}. To distinguish the new notions from the original one defined above, we sometimes call the new notion the  \emph{asymptotic profiles}. 

For example, let us consider the following  scalar quasilinear wave equation \fm{g^{\alpha\beta}(u)\partial_\alpha\partial_\beta u=0,\qquad\text{in }\R^{1+3},}
along with initial data \eqref{init}.
For sufficiently small $\eps>0$, this equation does admit a global solution $u$, but in general, the limit \eqref{def1a1} does not exist, not even if we replace the straight line with the characteristic. This is because a global solution to the wave equation above has a pointwise decay $\eps t^{-1+C\eps}$, and the $C\eps$ in the power of $t$ cannot be improved in general. See, e.g.,  \cite{1lind}. In this case, we construct  an asymptotic profile, not by taking the limit \eqref{def1a1}, but by solving a certain system of asymptotic equations. For example, we can take our asymptotic profile as a solution to the H\"{o}rmander's asymptotic equation \eqref{asypde11}, or to the geometric reduced system introduced by the author \cite{1yu2020}. 

In general, given a system of nonlinear wave equations \eqref{nlw} which admits a global solution for any given initial data \eqref{init}, we are  interested in finding a good notion of asymptotic profile. One could then ask the following two types of questions related to this asymptotic profile.
\begin{enumerate}[1.]
\item Given a global solution to \eqref{nlw}, can we find a corresponding asymptotic profile in the sense defined above? If two solutions correspond to the same asymptotic profile, are these two solutions the same?
\item Given an asymptotic profile, can we construct a global solution to \eqref{nlw}  which matches the asymptotic profile at the infinite time?
\end{enumerate}
In fact, these are the two main problems studied in (modified) scattering theory. The first one is called \emph{asymptotic completeness} and the second one is called \emph{existence of (modified) wave operators}. We refer our readers to \cite{1lindschl,1yu2020,yu2021asymptotic,1dengpusa,1dafermos2013scattering} for some work on (modified) scattering theory for nonlinear wave equations.}

\rm\bigskip

It is now natural for us to study the uniqueness properties related to the radiation fields. Suppose that two solutions to \eqref{swe} with (possibly different) $C_c^\infty$ initial data  and that their radiation fields are equal everywhere. We are now interested in whether these two solutions are also equal. If \eqref{swe} is replaced by $\Box u=0$, then the answer is yes; see Theorem 6.2.2 in \cite{1horm}. More interestingly, we can ask what happens if we only assume that the radiation fields of two solutions are equal whenever $q>R$ for some $R\in\R$. This would be the main question we  study in this paper, in the context of the equations \eqref{swe}.

\subsection{Unique continuation}\label{s1ucp}
The uniqueness results proved in this paper are usually referred to as \emph{unique continuation}. We refer to Tataru \cite{tata} for a survey on this topic. Generally, in unique continuation, we ask the following question:
\que{\label{ucp}\rm Let $P=\sum_{|\alpha|\leq m}c_\alpha(x)\partial^\alpha$ be an $m$-th order linear partial differential operator, and let $A$ and $B$ be two regions with $A\subset B$. Suppose that $u$ is a solution to $Pu=0$ in $B$ and that $u=0$ in $A$. Does it follow that $u=0$ in $B$?}\rm\\

For example, if $P$ is the Laplacian $\Delta$, $A$ is a nonempty open set and $B$ is a connected open set, then we have an affirmative answer to Question \ref{ucp}.

One is interested in this type of question especially when it is related to an ill-posed Cauchy problem. In many cases, a  Cauchy problem may not admit a solution for some initial or boundary data, but one could still expect a uniqueness result whenever a solution exists.

\subsubsection{Unique continuation  across a surface}\label{sec1.1.1}

There have been several (local) unique continuation results proved in the case when $A$ is given by a level set. To state these results, we  reformulate Question \ref{ucp} as follows.

\que{\label{ucp2} \rm Let $P=\sum_{|\alpha|\leq m}c_\alpha(x)(\partial/i)^\alpha$ be an $m$-th order linear partial differential operator. Let $h$ be a function and $S$ be a level set of $h$. Fix $x_0\in S$. Suppose that $u$ is a smooth ($C^\infty$) solution to $Pu=0$ in some neighborhood $V$ of $x_0$, and that $u=0$ in $V\cap\{h>h(x_0)\}$. Does it follow that $u=0$ near $x_0$?}\rm\\

If the coefficients of $P$ are all real analytic, and if $S$ is noncharacteristic with respect to $P$ at $x_0$ (i.e.\ $p(x_0,\nabla h(x_0))\neq 0$ where $p=p(x,\xi)$ is the principal symbol of $P$), then we  have an affirmative answer to Question \ref{ucp2}. This is   the Holmgren's theorem; see \cite{holmgren1901systeme,john2,john}. We also remark that this theorem is related to the Cauchy–Kowalevski theorem.

If the coefficients of $P$ are merely smooth, in general, we do not have unique continuation for non-characteristic surfaces. See \cite{cohen,horm2,alinbaou} for some counterexamples. To guarantee unique continuation, we need  additional assumptions.
 
\defn{\rm Suppose that $P$ is a linear differential operator defined in $\R^d$. Let $p$ be the principal symbol of $P$. That is, $p(x,\xi)=\sum_{|\alpha|=n}c_\alpha(x)\xi^\alpha$. We say that $P$ is \emph{principally normal} in an open set $X\subset\R^d$ if for any compact subset $K$ of $X$, we have
\fm{|\{\overline{p},p\}(x,\xi) |\lesssim_K |p(x,\xi)||\xi|^{m-1},\quad\forall (x,\xi)\in K\times\R^d.} 
Here $\{p,q\}$ is the Poisson bracket defined by $\{p,q\}=\partial_\xi p\cdot\partial_xq-\partial_x p\cdot \partial_\xi q$.}
\rmk{\rm In this paper we will only study unique continuation for operators with real principal symbols. From the definition, all such operators are principally normal.}\rm

\defn{\label{s1defnpsp}\rm Let $P$ be a principally normal operator defined in $\R^d$ whose principal symbol is $p$. Fix $x_0\in S$ where $S$ is a level set of a $C^2$ function $h$. Suppose that $\nabla h(x_0)\neq 0$. Then,  $S$ is \emph{strongly pseudoconvex} at $x_0$ with respect to $P$, if we have
\fm{{\rm Re}\{\overline{p},\{p,h\}\}(x_0,\xi)>0,\quad \text{if }\xi\in\R^d\setminus 0,\ p(x_0,\xi)=\{p,h\}(x_0,\xi)=0;}
\fm{&\{\overline{p}(x,\xi-i\tau\nabla h(x)),p(x,\xi+i\tau\nabla h(x))\}/2i\tau>0,\\ &\hspace{1cm}\text{if }x=x_0,\ \xi\in\R^d,\ \tau>0,\  p(x,\xi+i\tau\nabla h(x))=\{p(x,\xi+i\tau\nabla h(x)),h\}=0.}}
\rmk{\rm Intuitively, strong pseudoconvexity means that all characteristic curves tangent to the surface $S$ at $x_0$ must bend towards the region $\{h>h(x_0)\}$ where the solution is assumed to be vanishing.}\rm\\

If the operator $P$ is principally normal with $C^2$ coefficients and if the surface $S$ is strongly pseudoconvex, then we have an affirmative answer to Question \ref{ucp2}. This is the H\"{o}rmander's theorem, and we refer our readers to Chapter 28 of H\"{o}rmander \cite{hormIV}. We also remark that the necessity of strong pseudoconvexity for unique continuation is suggested by a counterexample constructed by Alinhac and Baouendi \cite{alinbaou}.

In fact, there is an intermediate case between the real analytic case and the smooth case. We decompose the space $\R^d$ as $\R^{d_1}\times\R^{d_2}$, and write $x\in\R^d$ as $x=(x_1,x_2)=\R^{d_1}\times\R^{d_2}$.  Assume that the coefficients of $P$ are partially analytic,  i.e.\ the coefficients $c_\alpha$'s are real analytic with respect $x_1$ and $C^1$ with respect to $x_2$. Then, under some suitable pseudoconvexity conditions, Tataru \cite{tata2,MR1697040}, H\"{o}rmander \cite{MR2033496} and Robbiano-Zuily \cite{robbzuil} proved that we still have an affirmative answer to Question \ref{ucp2} in this case. 

\subsubsection{Unique continuation for the wave equations}
In Section \ref{sec1.1.1}, when the coefficients of $P$ are partially analytic, we mentioned that there is a local unique continuation result proved in \cite{tata2,MR1697040,MR2033496,robbzuil}. In the wave equation setting, we formulate their results as follows. Let us consider the linear wave-type equation:
\fm{Pu:=\Box_gu+V(u)+Wu=0\qquad \text{in }\R^{1+d}_{t,x}.}
Here $g$ is a given Lorentzian metric, $\Box_g$ is the corresponding Laplace-Beltrami operator, $V$ is a vector field, and $W$ is a potential function. Suppose that the coefficients of $P$ are all smooth, and real analytic with respect to $t$. Then, we have local unique continuation  across any timelike surface.

We note that all the results in Section \ref{sec1.1.1} are local. That is, they only hold in a neighborhood of a certain point. There are in fact several  nonlocal uniqueness results in the wave equation setting. Some of them have the same form as Question \ref{ucp}. For example,  Ionescu and Klainerman \cite{ioneklai} proved a uniqueness theorem  for the wave equations across  bifurcate and characteristic surfaces $\{|x|=|t|+1\}\subset\R^{1+d}_{t,x}$; Whitman and P.\ Yu \cite{whityu} showed a converse theorem of the classical Huygens principle for free wave equations. There are also results which connect uniqueness with the   decays of the solutions at infinity.  For example,   Alexakis and Shao \cite{alexshao} proved that a solution to $\Box u+Vu=0$ must vanish if there is  no incoming and no outgoing radiation on specific halves of past and future null infinities;  Alexakis, Schlue and Shao \cite{alexschlshao}  proved   various uniqueness results from null infinity, for linear waves on asymptotically flat space-times; Duyckaerts, Kenig and  Merle \cite{MR4289254} (also see \cite{MR2966655,MR4163362} for some related work) proved that a radial nonradiative solution (i.e.\ a solution with asymptotically vanishing energy outside $|x|=|t|$  as $t\to\pm\infty$) to the energy-critical focusing wave equation with $C_c^\infty$ data vanishes everywhere. Since the radiation fields describe the asymptotic behavior of solutions to wave equations, we can see that the main question studied in this paper (stated right before Section \ref{s1ucp}) is closely related to the examples above.

\subsection{The main theorems}
We now state the main theorem of this paper.

\thm{\label{mthm0} Suppose that $u$ and $\wt{u}$ are two smooth global solutions to \eqref{swe} for $t\geq 0$ with $C_c^\infty$ initial data. Assume that the radiation fields $F_0$ and $\wt{F}_0$ of $u$ and $\wt{u}$, respectively, exist in the sense of Definition \ref{defrad}. Suppose that  $(u,u_t)|_{t=0}(x)=(\wt{u},\wt{u}_t)|_{t=0}(x)$ whenever $|x|\geq R_1$, and that $F_0(q,\omega)=\wt{F}_0(q,\omega)$ whenever $q>R_2$. Here  $R_1,R_2$ are real constants such that $R_1>0$ and $R_2<R_1$.

Our conclusion is that
\begin{enumerate}[\rm (i)]
\item If $|R_2|<R_1$, then $u=\wt{u}$ whenever $t>0$ and $|x|^2-(t+R_2)^2\geq R_1^2-R_2^2$. 
\item If $R_2\leq -R_1$, then $u=\wt{u}$ everywhere.
\end{enumerate}}\rm\bigskip

We quickly remark that part (i) of Theorem \ref{mthm0} can be improved in the linear case (i.e.\ $Q\equiv 0$ in \eqref{swe}). We will discuss this later in Remark \ref{rmkk1.4}.

To prove Theorem \ref{mthm0}, we first notice that part (ii) follows from part (i) by sending $R_2\downarrow -R_1$. To prove part (i), we  study the linear  PDE which the difference $u-\wt{u}$ should satisfy. For an $\R^N$-valued function $\phi=(\phi^I)_{I=1}^N$, we define $\wt{\Box}\phi=((\wt{\Box}\phi)^I)_{I=1}^N$ by
\eq{\label{lpddefn}(\wt{\Box}\phi)^I&:=\Box \phi^I+\sum_{J=1}^N\sum_{\alpha=0}^3V_J^{I,\alpha} \partial_\alpha \phi^J+\sum_{J=1}^NW_J^I\phi^J,\qquad I=1,2,\dots,N.}
Here $V_{J}^{I,\alpha}=V_{J}^{I,\alpha}(t,x)$ and $W_J^I=W_J^I(t,x)$ are given continuous functions.  

Theorem \ref{mthm0} is now a corollary of the following theorem.
\thm{\label{mthm} Fix two constants $R_1,R_2\in\R$ such that $R_1>0$ and $|R_2|<R_1$. Define an open set
\eq{\label{dddefn} \D=\D_{R_1,R_2}:=\{(t,x)\in\R^{1+3}:\ t> 0,\ (r-t-R_2)(r+t+R_2)>R_1^2-R_2^2\}.}
For a fixed constant $0<\gamma\leq 1$, we suppose that the continuous functions $V_{J}^{I,\alpha}$'s and $W^{I}_J$'s for $I,J=1,\dots,N$ and $\alpha=0,1,2,3$ are defined in $\D$ and that they satisfy the following pointwise estimates:
\eq{\label{mthma1}\sum_{I,J=1}^N\sum_{\alpha=0}^3|V^{I,\alpha}_J|\leq C_1 (1+t+r)^{-1}\qquad\forall(t,x)\in\D,}
and
\eq{\label{mthma2}\sum_{I,J=1}^N|\sum_{\alpha=0}^3V^{I,\alpha}_J\widehat{\omega}_\alpha|+\sum_{I,J=1}^N|W^{I}_J|\leq  C_1(1+t+r)^{-1-\gamma}\qquad\forall(t,x)\in\D.}
Here $\widehat{\omega}=(-1,\omega)=(-1,x/r)$.

In addition, we suppose that $\phi=(\phi^I)\in C^2(\D;\R^N)$ is a solution to $\wt{\Box}\phi=0$.
Assume that $\phi=0$ for all $r-t\geq R_1$, and that for a fixed constant $\gamma'>0$, it satisfies the pointwise estimates: 
\eq{\label{mthma3}|\phi|+|\partial\phi|\leq C_2(1+t+r)^{-3/2-\gamma'}\qquad \forall(t,x)\in\D,}
and
\eq{\label{mthma4}|(\partial_t+\partial_r)\phi|\leq C_2(1+t+r)^{-5/2-\gamma'}\qquad \forall(t,x)\in\D.}

Our conclusion is that  $\phi\equiv 0$ in $\D$.}\rm\bigskip

\begin{figure}
    \centering

\begin{tikzpicture}
\begin{axis}[
    axis equal image,
    axis lines = middle,
    xlabel = \(r\),
    ylabel = {\(t\)},
    yticklabels={,,},
    xticklabels={,,}
]

\addplot [
    name path = A,
    domain=3:15, 
    samples=100, 
    color=red
]
{(x^2-8)^(1/2)-1};
\addplot [
    name path= D,
    domain=1:15, 
    samples=100, 
    color=blue,
    dashed
    ]
    {x-1};
\addplot [
    name path = B,
    domain=3:15, 
    samples=100, 
    color=blue,
    ]
    {x-3};
\addplot [
    name path = C1,
    domain=0:3, 
    samples=100, 
    color=blue,
    dashed
    ]
    {3-x};
\addplot [red!25] fill between [of = A and B, soft clip={domain=3:15}];
\addplot[
    name path = E,
    domain=3:15, 
    samples=100, 
    color=black
    ]
    {0};
\addplot [blue!25] fill between [of = B and E, soft clip={domain=3:15}];
\addplot [green!25] fill between [of = A and D, soft clip={domain=3:15}];
\addplot [green!25] fill between [of = C1 and D, soft clip={domain=2:3}];

\addplot[color=blue!25] (10,2) node[color=black] {\tiny $A$};
\addplot[color=red!25] (10,7.9) node[color=black] {\tiny $B$};
\addplot[color=green!25] (3,1.2) node[color=black] {\tiny $C$};

\addplot[color=white,domain=15:16] {-1};
\addplot[color=white,domain=-1:-0.5] {-1};
\addplot[color=black] (0,0) node[above left,color=black] {$0$};
\addplot (3,0) node[below] {\tiny $R_1$};
\addplot (0,3) node[below] {\tiny $R_1$};
\addplot (1,0) node[below] {\tiny $R_2$};
\end{axis}
\end{tikzpicture}

\caption{Assuming $|R_2|<R_1$, we define the blue region, the red region and the green region  by $A=\{r-t\geq R_1\}$, $B=\{ r^2-(t+R_2)^2\geq R_1^2-R_2^2\}\setminus A$,  $C=\{r-t\leq R_2,\ r+t\geq R_1\}\setminus(A\cup B)$, respectively.}
\label{fig}
\end{figure}
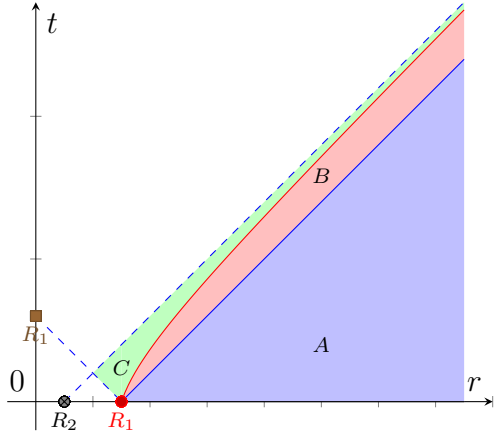

In Figure \ref{fig}, we depict the regions  in the main theorems.  The union of the red and blue regions $A\cup B$ denotes both the set where  $u=\wt{u}$ in part (i) of Theorem \ref{mthm0}, and the region $\D$ in Theorem \ref{mthm}.

Here are several remarks about Theorem \ref{mthm0} and \ref{mthm}.

\rmkk{\rm In Theorem \ref{mthm0}, we do not have any smallness assumption; the small parameter $\eps$ is not involved at all. As a result, this theorem is also applicable in the large data case. However, for large data, a global solution does not necessarily exist; even if it exists, it does not necessarily admit a radiation field. Thus, we must assume that there are two global solutions and that these global solutions admit radiation fields in the statement of Theorem~\ref{mthm0}.}

\rmkk{\label{rmkk1.4}\rm
Both our main theorems involve a hyperboloid $r^2-(t+R_2)^2=R_1^2-R_2^2$. This hyperboloid arises as follows. In our proofs, we apply the H\"{o}rmander's theorem to a family of  strongly pseudoconvex surfaces $(r+\nu)^2-(t+R_2)^2=R_1^2-R_2^2$ with $\nu>0$. See Example \ref{exmhyp}. By taking the limit as  $\nu\downarrow 0$, we obtain the hyperboloid above. However, this hyperboloid itself  is not  strongly pseudoconvex. In other words, the exterior of this hyperboloid is the largest region we 
expect in the uniqueness result  if we only apply the H\"{o}rmander's theorem. 

In some special cases, the conclusions of our main theorems can be improved.  In Theorem \ref{mthm}, if all $V_*^{**},W_{*}^*\equiv 0$ and if $|R_2|<R_1$, we  can replace $\D$ with a larger region \fm{\{(t,x)\in\R^{1+3}:\ t>0,\ r+t\geq R_1,\ r-t\geq R_2\}}which is the union  $A\cup B\cup C$ in Figure \ref{fig}.  Similarly, we have a better result in Theorem \ref{mthm0} if we assume that all $Q^*\equiv 0$. This uniqueness result for the linear wave equation, to the author's knowledge, seems a new one.  We shall prove this result at the end of Section \ref{sec5}. However, in the proof we need to apply the Holmgren's theorem, so unfortunately it is unclear whether the conclusions can be improved in general.}

\rmkk{\rm Theorem \ref{mthm} would fail if we do not assume that $\phi$ is vanishing for $r-t\geq R_1$. In fact,  $1/r$ is a solution to the linear wave equation in $\R^{1+3}\setminus \{x=0\}$, and so is any derivative of $1/r$. Now, for each $k\geq 1$, we have $|\partial^k (1/r)|\lesssim r^{-1-k}$, so its radiation field is zero. However,  $\partial^k (1/r)$ is nonvanishing everywhere.}

\rmkk{\rm  In Theorem \ref{mthm0}, we cannot obtain a uniqueness result if we only have information about where the two radiation fields are equal. In other words, the assumption that $(u,u_t)|_{t=0}=(\wt{u},\wt{u}_t)|_{t=0}$ whenever $|x|\geq R_1$  is necessary for our uniqueness result.

To see this, we  construct a family of solutions to the linear wave equation $\Box u=0$, such that all of their radiation fields vanish for $q\geq 0$ and that the union of their supports is equal to $\R^{1+3}$. Suppose that $u$ is a nonzero  solution to $\Box u=0$ such that the support of $(u,u_t)|_{t=0}$ is $\{r\leq 1\}$. Because of the finite speed of propagation, the support of $u$ is $\{r-|t|\leq 1\}$. Set $u_L(t,x):=u(t-L,x)$, and we get a family of solutions to the linear wave equation. For all $L>1$,  the radiation field $F_{L,0}(q,\omega)$ of $u_L$ is equal to $0$ whenever $q\geq 0$, and the support of $u_L$ is $\{r-|t-L|\leq 1\}$. It is clear that $\bigcup_{L>1}\{r-|t-L|\leq 1\}=\R^{1+3}$. }
\rm

\subsection{Idea of the proofs}
We first discuss how to prove Theorem \ref{mthm}. Define a function \fm{f_\delta(t,x)=r-t-\delta^{-1}t^{-\delta},\qquad 0<\delta<1/2,\ t>0.} We choose this function because its level set is strongly pseudoconvex with respect to $\Box$ for sufficiently large time $t$, and because its level set approaches the characteristic level set $\{r-t=C\}$ as $t\to\infty$. In Section \ref{s3}, we prove the Carleman estimates, i.e.\ Proposition \ref{prop3.1}, which are related to this function $f_\delta$. The proof is similar to that of classical Carleman estimates, but we need to take the decay rates of all the functions involved into account. In Section \ref{sec4}, we show a unique continuation result related to $f_\delta$ by applying the Carleman estimates. Note that Proposition \ref{prop3.1} only works for functions which are compactly supported in the spacetime, so we  introduce some cutoff functions in our proof. The main result in this section, Proposition \ref{prop4.1}, states that a solution $\phi$ in Theorem \ref{mthm} must vanish in a region of the form $\{t>T_0,\ f_\delta>R''\}$ for some well chosen constants $T_0,R''$, under reasonable assumptions. Then, in Section \ref{sec5}, we prove Theorem \ref{mthm} by applying Proposition \ref{prop4.1} and  Corollary \ref{cors5carest} which is the H\"{o}rmander's theorem.

After proving Theorem \ref{mthm}, we notice that Theorem \ref{mthm0} follows from Theorem \ref{mthm}. The proof relies on both the null condition \eqref{swea2} and the estimate \eqref{def1a2} in Definition \ref{defrad}. See Section \ref{sec5.3}.  We also explain how we prove the result stated in Remark \ref{rmkk1.4} by applying the Holmgren's theorem. It turns out that the proof of this part is very similar to that of Theorem \ref{mthm0}.

\subsection{Acknowledgement}
The author would like to thank his advisor, Daniel Tataru, for suggesting this problem and for many helpful discussions. The author would like to thank Shi-Zhuo Looi, for pointing out Remark \ref{mrmk1.10}.  The author is also grateful
to the anonymous reviewer for several valuable comments and suggestions on this paper.

This research was partially supported by a James H.\  Simons Fellowship, by the NSF grant DMS-1800294, and by the Deutsche Forschungsgemeinschaft (DFG, German Research Foundation) under
Germany's Excellence Strategy - GZ 2047/1, Projekt-ID 390685813. 

\section{Preliminaries}\label{sec2}
\subsection{Notations} We use $C$ to denote universal positive constants. We write $A\lesssim B$ or $A=O(B)$ if $|A|\leq CB$ for some  $C>0$. We write $A\sim B$ if $A\lesssim B$ and $B\lesssim A$. We  use $C_{v}$ or $\lesssim_v$ if we want to emphasize that the constant depends on a parameter $v$. The values of all constants in this paper may vary from line to line. Moreover, we write $A\ll 1$  if $0<A<1$ is a sufficiently small number, and we write $A\gg 1$  if $A>1$ is a sufficiently large number. Similarly, we use $\ll_v$ and $\gg_v$ to emphasize the dependence on a parameter $v$.

Unless specified otherwise, we always assume that the Latin indices $i,j,l$ take values in $\{1,2,3\}$ and the Greek indices $\alpha,\beta$ take values in $\{0,1,2,3\}$.  We define $\partial_\alpha$ by $\partial_0=\partial_t$ and $\partial_i=\partial_{x_i}$ for $i=1,2,3$. We also define $\partial^\alpha=m^{\alpha\beta}\partial_\beta$ with $(m^{\alpha\beta})=(m_{\alpha\beta})=\text{diag}(-1,1,1,1)$ is the Minkowski metric. Here we use the Einstein summation convention. In other words, we have $\partial^0=-\partial_0$ and $\partial^i=\partial_i$.

\subsection{Commuting vector fields}\label{sec2.1}

We denote by $Z$ any of the following vector fields:
\begin{equation}\label{vf} \partial_\alpha,\ \alpha=0,1,2,3;\ S=t\partial t+r\partial_r;\ \Omega_{ij}=x_i\partial_j-x_j\partial_i,\ 1\leq i<j\leq 3;\ \Omega_{0i}=x_i\partial_t+t\partial_i,\ i=1,2,3.\end{equation}We write these vector fields as $Z_1,Z_2,\dots,Z_{11}$, respectively. For any multiindex $I=(i_1,\dots,i_m)$ with length $m=|I|$ such that $1\leq i_*\leq 11$, we set $Z^I=Z_{i_1}Z_{i_2}\cdots Z_{i_m}$. Then we have the Leibniz's rule
\begin{equation}Z^I(fg)=\sum_{|J|+|K|=|I|}C^I_{JK}Z^JfZ^Kg,\hspace{1cm}\text{where $C_{JK}^I$ are constants.}\end{equation}
\rmkno{\rm In this paper, the superscript $I$ can denote either a multiindex (e.g.\ $Z^I$) or an index for an $\R^N$-valued function (e.g.\ $u^I$).  }\rm\bigskip

We have the following commutation properties.
\begin{equation}[S,\square]=-2\square,\hspace{1cm}[Z,\square]=0\text{ for other $Z$};\end{equation}
\begin{equation}[Z_1,Z_2]=\sum_{|I|=1} C_{Z_1,Z_2,I}Z^I,\hspace{1cm}\text{where $C_{Z_1,Z_2,I}$ are constants};\end{equation}
\begin{equation}\label{comf3}[Z,\partial_\alpha]=\sum_\beta C_{Z,\alpha\beta}\partial_\beta,\hspace{2em} \text{where $C_{Z,\alpha\beta}$ are constants}.\end{equation}

In this paper, we will need the following pointwise bounds.

\begin{lem}\label{l2.1}
For any function $\phi$, we have 
\begin{equation}|\partial^k\phi|\leq C\lra{t-r}^{-k}\sum_{|I|\leq k}|Z^I\phi|,\hspace{1cm}\forall k\geq 1,\end{equation}
and
\begin{equation}\label{l2.1c1}|(\partial_t+\partial_r)\phi|+|(\partial_i-\omega_i\partial_r)\phi|\leq C\lra{t+r}^{-1}\sum_{|I|=1}|Z^I\phi|.\end{equation}
Here $\lra{s}=\sqrt{1+|s|^2}$ is the Japanese bracket.
\end{lem}\rm

The proofs of  these results are standard. We refer  to, for example, \cite{1sogg,1horm} for their proofs.

\subsection{The null condition}
Let $Q=(Q^I)$ be a smooth $\R^N$-valued  function defined in  $\R^{N}\times\R^{4N}$, and suppose that the nonlinear terms $Q^I$'s satisfy the null condition. That is, we have \eqref{swea1} and \eqref{swea2}. Then we have the following useful estimates.
\prop{\label{nullest} Fix an integer $M>0$ and a time $T_0>0$. Suppose $\phi$ is an $\R^N$-valued $C^M$ function of $(t,x)$   such that 
\eq{\label{nullestasu}\sum_{|J|\leq M}|Z^J\phi|\leq 1,\qquad\qquad\text{for }t\geq T_0.}
Then, whenever $r/ t\in[1/2,2]$ and $t\geq T_0$, we have
\eq{\label{nullestcon}\sum_{|J|\leq M-1}|Z^JQ(\phi,\partial\phi)|\lesssim_{M} \lra{t}^{-1}\lra{r-t}^{-1}(\sum_{1\leq |J|\leq M}|Z^J\phi|)^2+(\sum_{|J|\leq M}|Z^J\phi|)^3.}}
\begin{proof}(Sketch). It suffices to prove \eqref{nullestcon} when there is no remainder term in \eqref{swea1} for each $I$. In other words, we can assume that
\fm{Q^I(u,v)&=\sum_{I,J=1}^N\sum_{\alpha,\beta=0}^3A^{\alpha\beta}_{I,JK}v_\alpha^J v_\beta^K,\qquad\text{for all } (u,v)\in\mcl{O}\subset\R^N\times\R^{4N}.}
We can prove this by the chain rule, the Taylor expansion of $Q$ and \eqref{nullestasu}. It thus suffices to prove that for fixed $I,J,K$, we have
\eq{\label{nullestpf1}\sum_{|L|\leq {M-1}}|Z^L(A_{I,JK}^{\alpha\beta}\partial_\alpha \phi\partial_\beta\phi)|\lesssim_M \lra{t}^{-1}\lra{r-t}^{-1}(\sum_{1\leq |J|\leq M}|Z^J\phi|)^2.}
Here we use the Einstein summation convention. This estimate follows from Lemma 6.6.5 in \cite{1horm}, Lemma II.5.4 in \cite{1sogg}, and Lemma \ref{l2.1}. We also refer to Lemma 2.3 and Lemma 2.4 in \cite{yu2021asymptotic}.
\end{proof}\rm

\subsection{Global existence results for \eqref{swe} with small $C_c^\infty$ data}\label{sec2.2}
In this subsection, we present several results about the global solutions to the Cauchy problem \eqref{swe} along with the initial data \eqref{init}. In fact we will not use the results in this subsection to prove our main theorems. Our goal, however, is to verify that Definition \ref{defrad}  makes sense. That is, a solution to \eqref{swe} does admit a radiation field in the sense of Definition \ref{defrad}, at least when the initial data are small and localized.

The first result is the global existence result for semilinear wave equations satisfying the null condition. In addition to \cite{klai,klai2,chri}, we also refer our readers to Section 6.6 of \cite{1horm}.

\prop[Klainerman \cite{klai,klai2}, Christodoulou \cite{chri}]{\label{gloexi}For sufficiently small $\eps\ll1$, the Cauchy problem \eqref{swe} and \eqref{init} has a smooth global solution $u$ for $t\geq 0$. Moreover, for each integer $M>0$, as long as $\eps\ll_M1$, we have the pointwise bounds
\eq{\label{gloexif1}\sum_{|J|\leq M}|Z^Ju(t,x)|\lesssim_{M}\eps\lra{t+r}^{-1}.}
Here $Z$'s are the commuting vector fields introduced in Section \ref{sec2.1}.}\rm\bigskip

Next, we present a global solution constructed in Proposition \ref{gloexi} admits a radiation field as defined in Definition \ref{defrad}. 

\prop{\label{radfie} Fix an integer $M>0$, and let $u$ be a smooth solution constructed in Proposition \ref{gloexi} satisfying the pointwise bounds \eqref{gloexif1}. Then, $u$ has a radiation field $F_0$ defined by 
\eq{\label{radfie1}F_0(q,\omega)=\lim_{r\to\infty} ru(r-q,r\omega),}
such that for $\eps\ll_M1$, we have
\eq{\label{radfie2}\sum_{|J|\leq M}|Z^J(u-r^{-1}F_0(r-t,\omega))|\lesssim_M \eps \lra{t}^{-2}\lra{r-t},\qquad r> t/2\gtrsim 1.}}
\begin{proof} (Sketch). Fix a multiindex $L$. Then, by applying $Z^L$ to \eqref{swe}, we can show that
\fm{\Box Z^Lu^I=Z^LQ^I(u,\partial u)+\sum_{|L'|<|L|}C_{L',L}Z^{L'}Q^I(u,\partial u).}
Then, for $\eps\ll_M1$, by  Proposition \ref{nullest} we have
\fm{\sum_{|L|\leq M}|\Box Z^Lu|\lesssim_M \eps^2 \lra{t+r}^{-3},\qquad\qquad t\gg1.}
Recall that \fm{\Box\phi&=-r^{-1}(\partial_t+\partial_r)(\partial_t-\partial_r)(r\phi)+r^{-2}\sum_{1\leq i<j\leq 3}\Omega_{ij}\Omega_{ij}\phi.} By setting $\phi=Z^Lu$ and applying \eqref{gloexif1}, we have
\fm{\sum_{|L|\leq M}|(\partial_t+\partial_r)[(\partial_t-\partial_r)(rZ^Lu)](t,x)|\lesssim \eps \lra{t+r}^{-2}.}
For each fixed $q_0$ and $\omega\in\mathbb{S}^2$, we can show that whenever $1\lesssim_{q_0} R<R'$,
\fm{&\sum_{|L|\leq M}|[(\partial_t-\partial_r)(rZ^Lu)](R'-q_0,R'\omega)-[(\partial_t-\partial_r)(rZ^Lu)](R-q_0,R\omega)|\\&\lesssim \int_{R}^{R'}\eps\lra{r}^{-2}\ dr\lesssim \eps R^{-1}.}
As a result, for each multiindex $L$ with $|L|\leq M$, the limit
\fm{A_{0,L}(q_0,\omega):=\lim_{r\to\infty} [(\partial_t-\partial_r)(rZ^Lu)](r-q_0,r\omega)}
exists. Thus, we can define 
\fm{F_0(q_0,\omega):=-\frac{1}{2}\int_{q_0}^\infty A_{0,0}(q,\omega)\ dq.}
The  proofs of \eqref{radfie1} and \eqref{radfie2} are standard. We refer to Section 5 and 7 in \cite{yu2021asymptotic} for similar proofs in a quasilinear wave equation setting.
\end{proof}
\rm

\subsection{The local unique continuation property}\label{sec5.1}
In this subsection, we recall the definition of strong pseudoconvexity and the Carleman estimates  for $P=\wt{\Box}$ where $\wt{\Box}$ is defined by \eqref{lpddefn}.   In particular, we state the corresponding H\"{o}rmander theorem (Corollary \ref{cors5carest}) which will be used in the proof of our main theorems; see Section \ref{sec5}. 

Set $p(\xi)=-m^{\alpha\beta}\xi_\alpha\xi_\beta$ for each $\xi=(\xi_\alpha)_{\alpha=0,1,2,3}\in\R^{1+3}$, so we have $p(D)=p(\nabla_{t,x}/i)=\Box$. We first recall the definition of strong pseudoconvexity.

\defn{\rm Fix an open subset $\mcl{O}$ in $\R^{1+3}$, and fix a $C^2$ function $\wt{f}(t,x)$ defined in  $\mcl{O}$. Suppose that $\nabla_{t,x}\wt{f}$ is nowhere vanishing. We say that the level sets of $\wt{f}$ are \emph{strongly pseudoconvex} with respect to $\Box$ in $\mcl{O}$, if the following two estimates hold for all $(t,x)\in\mcl{O}$: we have \eq{\label{spdefnf1}\{p,\{p,\wt{f}\}\}(t,x,\xi)>0}
whenever $\xi\in\R^{1+3}\setminus\{0\}$ and $p(\xi)=\{p,\wt{f}\}(t,x,\xi)=0$;  for each $\tau>0$ and $p_{\wt{f}}(t,x,\xi):=p(\xi+i\tau\nabla_{t,x}\wt{f}(t,x))$, we have
\eq{\label{spdefnf2}\frac{1}{i\tau}\{\overline{p_{\wt{f}}},p_{\wt{f}}\}(t,x,\xi)>0}whenever $\xi\in\R^{1+3}\setminus\{0\}$ and
$p_{\wt{f}}(t,x,\xi)=\{p_{\wt{f}},\wt{f}\}(t,x,\xi)=0$.}
\rmkno{\rm Since $p(\xi)$ is $\R$-valued, we can prove that \eqref{spdefnf1} implies \eqref{spdefnf2}. Thus, to prove strong pseudoconvexity, we  only need to check \eqref{spdefnf1}.}\rm\bigskip

Strong pseudoconvexity implies the Carleman estimates, as shown in the next proposition. Note that $p(\xi)$ is the principal symbol of $\wt{\Box}$ where $\wt{\Box}$ is defined by \eqref{lpddefn}.

\prop[Carleman estimates]{\label{s5carest} Suppose that the level sets of $\wt{f}$ are strongly pseudoconvex with respect to $\Box$ in an open subset $\mcl{O}$. Fix $(t_0,x_0)\in\mcl{O}$. Then, there exist constants $\tau_0,\lambda>0$ and an open set $\mcl{N}\subset\mcl{O}$ containing $(t_0,x_0)$, such that 
\eq{\label{s5carestf1} \tau\norm{\exp(\tau\exp(\lambda \wt{f}))\psi}_{H^1_\tau(\mcl{N})}^2\lesssim \norm{\exp(\tau\exp(\lambda \wt{f}))\wt{\Box} \psi}_{L^2(\mcl{N})}^2}
for each $\psi\in H^1_{c}(\mcl{N};\R^N)$ (i.e.\ $\psi$ is $\R^N$-valued and $\supp\psi$ is compact) and $\tau\geq\tau_0$. Here the $H^1_\tau$ norm is defined by \fm{\norm{g}_{H^1_\tau}:=\norm{(|D|^2+\tau^2)^{1/2}g}_{L^2}.}}
\begin{proof}We first note that the proposition holds with \eqref{s5carestf1} replaced by
\eq{\label{s5carestf2} \tau\norm{\exp(\tau\exp(\lambda \wt{f}))u}_{H^1_\tau(\mcl{N})}^2\lesssim \norm{\exp(\tau\exp(\lambda \wt{f}))\Box u}_{L^2(\mcl{N})}^2}
for each $u\in H^1_{c}(\mcl{N};\R)$. The proof of \eqref{s5carestf2} relies on  strong pseudoconvexity of $\wt{f}$. We skip its proof here and  refer our readers to Section 28.3 in H\"{o}rmander \cite{hormIV}, or Remark 7 and Theorem 8 in Tataru \cite{tata}.

Next, we write $\Phi:=\exp(\lambda\wt{f})$. By applying \eqref{s5carestf2} to $\psi^I$, we have
\fm{\tau^{1/2}\norm{\exp(\tau\Phi)\psi}_{H^1_\tau}&\lesssim\sum_I\norm{\exp(\tau\Phi)\Box (\psi^I)}_{L^2}\\&\lesssim \norm{\exp(\tau\Phi)\wt{\Box} \psi  }_{L^2}+\norm{\exp(\tau\Phi)(\sum_{J,\alpha}V_J^{I,\alpha}\partial_\alpha\psi^J+\sum_JW^{I}_J\psi^J)}_{L^2}\\
&\lesssim\norm{\exp(\tau\Phi)\wt{\Box} \psi  }_{L^2}+\norm{\exp(\tau\Phi)\partial\psi}_{L^2}+\norm{\exp(\tau\Phi)\psi}_{L^2}.}
Here we use the estimates $V_{*,*}^*,W^{*}_*=O(1)$. To continue, we note that
\fm{\norm{\exp(\tau\Phi)\partial\psi}_{L^2}\leq \norm{\exp(\tau\Phi)\psi}_{\dot{H}^1}+\norm{\tau\exp(\tau\Phi)(\partial\Phi)\psi}_{L^2}.}
Since $\partial\Phi=\lambda\exp(\lambda\wt{f})\cdot\partial\wt{f}=\lambda\Phi\partial\wt{f}=O_\lambda(1)$, we conclude that 
\fm{\tau^{1/2}\norm{\exp(\tau\Phi)\psi}_{H^1_\tau}&\leq C\norm{\exp(\tau\Phi)\wt{\Box} \psi  }_{L^2}+C_\lambda\norm{\exp(\tau\Phi)\psi}_{H^1_\tau}.}
Thus, by choosing $\tau_0\geq 4C_\lambda^2$, we have $\tau^{1/2}\geq 2C_\lambda$ whenever $\tau\geq \tau_0$. As a result, we have
\fm{\tau^{1/2}\norm{\exp(\tau\Phi)\psi}_{H^1_\tau}&\lesssim_\lambda\norm{\exp(\tau\Phi)\wt{\Box} \psi  }_{L^2},\qquad \tau\geq \tau_0.}
This is \eqref{s5carestf1}.
\end{proof}
\rm\bigskip

A corollary of the Carleman estimates is the following local unique continuation property. It is in fact the H\"{o}rmander's theorem stated after Definition \ref{s1defnpsp} in Section \ref{sec1.1.1} with $P=\wt{\Box}$. Again, we refer our readers to \cite{hormIV,tata}.
\cor[Local unique continuation property]{\label{cors5carest}Suppose that the level sets of $\wt{f}$ are strongly pseudoconvex with respect to $\Box$ in an open set $\mcl{O}$. Fix $(t_0,x_0)\in\mcl{O}$ and suppose $\phi$ is a solution to $\wt{\Box}\phi=0$ near $(t_0,x_0)$. Then, there exists an open set $\mcl{N}$ in $\R^{1+3}$  such that $(t_0,x_0)\in\mcl{N}\subset\mcl{O}$, such that, if  $\phi\equiv 0$ in $\mcl{N}\cap\{\wt{f}\geq \wt{f}(t_0,x_0)\}$, then  $\phi\equiv 0$ in $\mcl{N}$.}\rm
\bigskip

We end this subsection with some examples of strongly pseudoconvex surfaces. In this paper, we study those surfaces which are level sets of functions of the form $\wt{f}(t,x)=h(t,r)$. That is, $\wt{f}$ is spherically symmetric with respect to $x$.  In the next lemma, we present a sufficient condition for strong pseudoconvexity.

\lem{\label{lem5.3}Suppose $\wt{f}(t,x)=h(t,r)$ where $h$ is a $C^2$ function such that $\partial_rh\neq 0$ everywhere. Then, we have
\eq{\{p,\wt{f}\}=2\xi_0 \partial_th-2\sum_{j=1}^3\xi_j\omega_j\partial_rh,}
\eq{&\{p,\{p,\wt{f}\}\}\\&=4\xi_0^2(\partial_t^2h)-8\xi_0 \sum_{j=1}^3\xi_j\omega_j(\partial_t\partial_rh)+4(\sum_{j=1}^3\xi_j\omega_j)^2(\partial_r^2h)+4r^{-1}((\sum_{j=1}^3\xi_j^2-(\sum_{j=1}^3\xi_j\omega_j)^2)\partial_rh.}

If moreover $p(\xi)=\{p,\wt{f}\}(t,x,\xi)=0$ for some $\xi\in\R^{1+3}\setminus\{0\}$ and some $(t,x)$, then we have $\xi_0^2=\sum_{j=1}^3\xi_j^2>0$ and $\sum_{j=1}^3\xi_{j}\omega_j=(\partial_th/\partial_rh)\xi_0$. In this case, we have
\eq{\label{lemspcompf1}\{p,\{p,\wt{f}\}\}&=4\xi_0^2[\partial_t^2h-2 (\partial_th/\partial_rh)\partial_t\partial_rh+(\partial_th/\partial_rh)^2\partial_r^2h+r^{-1}(1-(\partial_th/\partial_rh)^2)\partial_rh].}

As a result, we have strong pseudoconvexity if we have  \eq{\label{lem5.3f1} \partial_t^2h-2 (\partial_th/\partial_rh)\partial_t\partial_rh+(\partial_th/\partial_rh)^2\partial_r^2h+r^{-1}(1-(\partial_th/\partial_rh)^2)\partial_rh>0}whenenver $|\partial_th|\leq |\partial_rh|$.}
\begin{proof}
Its proof is standard and thus omitted. We remark that, if $|\partial_th|>|\partial_rh|$, then $p(\xi)=\{p,\wt{f}\}(t,x,\xi)=0$ will not hold for any $(t,x,\xi)$. This is because $|\xi_0|=(\sum_j\xi_j^2)^{1/2}\geq|\sum_j\xi_j\omega_j|$ whenever $p(\xi)=0$.
\end{proof}
\rm\bigskip

Here are two  examples of $\wt{f}$ whose level sets are strongly pseudoconvex in $\mcl{O}$. Note that the strong pseudoconvexity proved  in Example \ref{exmhyp} will be used in Section \ref{sec5}; see \eqref{f511}.
\exm{\rm Let $\mcl{O}=\D$ and $h(t,r)=r-\nu t$ where $\nu\in\R\setminus\{\pm 1\}$ is a fixed constant. In fact, when $|\nu|>1$, we have $|\partial_th|>|\partial_rh|$; when $|\nu|<1$, we have 
\fm{\partial_t^2h-2 (\partial_th/\partial_rh)\partial_t\partial_rh+(\partial_th/\partial_rh)^2\partial_r^2h+r^{-1}(1-(\partial_th/\partial_rh)^2)\partial_rh=r^{-1}(1-\nu^2) >0.}
So the level sets of $h$ are strongly pseudoconvex in $\mcl{O}$.}
\exm{\label{exmhyp}\rm Fix $\nu>0$ and $\kappa\in\R$. Define $h(t,r)=(r+\nu)^2-(t+\kappa)^2$ and \fm{\mcl{O}=\{(t,x)\in\R^{1+3}:\ t>0,\ h>0\}.}
 Then,\fm{&\partial_t^2h-2 (\partial_th/\partial_rh)\partial_t\partial_rh+(\partial_th/\partial_rh)^2\partial_r^2h+r^{-1}(1-(\partial_th/\partial_rh)^2)\partial_rh\\
&=-2+2(\frac{2(t+\kappa)}{2(r+\nu)})^2+r^{-1}(1-(\frac{2(t+\kappa)}{2(r+\nu)})^2)\cdot 2(r+\nu)=2(r+\nu)^{-2} r^{-1}\nu h.}
 So the level sets of $h$ are strongly pseudoconvex in $\mcl{O}$.
}\rm

\section{The   Carleman estimates}\label{s3}

In this section, we present the Carleman estimates which are used to prove our main theorems. To state the main proposition, we first make  several additional definitions. Fix two fixed constants $R\in\R$ and $T_0>|R|$. We define
\eq{\label{regdef}\Omega_{T_0,R}:=\{(t,x)\in\R^{1+3}:\ t\geq T_0,\ |x|-t>R\}.}For a fixed constant $\delta\in(0,1/2)$, we set
\eq{\label{fdefn}f_\delta=f_\delta(t,x):=r-t-\delta^{-1}t^{-\delta},\qquad (t,x)\in\R_+\times\R^3,\ r=|x|.}
We usually omit the subscript and write $f$ instead of $f_\delta$.

In this section we seek to prove the following Carleman estimates.
\prop{\label{prop3.1}Fix $R,R',R''\in\R$ such that $R>R'>R''$. Also fix $T_s>|R''|+1$. For a fixed constant $0<\gamma\leq 1$, we suppose that we have continuous functions $V_{J}^{I,\alpha}$'s and $W^{I}_J$'s for $I,J=1,\dots,N$ and $\alpha=0,1,2,3$ defined in $\Omega_{T_s,R''}$, and that they satisfy the following pointwise estimates:
\eq{\label{prop3.1a1}\sum_{I,J=1}^N\sum_{\alpha=0}^3|V^{I,\alpha}_J|\leq C_1 (1+t+r)^{-1}\qquad\forall(t,x)\in\Omega_{T_s,R''},}
and
\eq{\label{prop3.1a2}\sum_{I,J=1}^N|\sum_{\alpha=0}^3V^{I,\alpha}_J\widehat{\omega}_\alpha|+\sum_{I,J=1}^N|W^{I}_J|\leq  C_1(1+t+r)^{-1-\gamma}\qquad\forall(t,x)\in\Omega_{T_s,R''}.}
Here $\widehat{\omega}=(-1,\omega)=(-1,x/r)$.

Then, for all sufficiently large time $T_0\gg_{T_s,R,R'',C_1,\gamma}1$, we have the following estimate. For all $\tau\gg 1$, $\delta\in(0,\gamma )$ and any function $\psi=(\psi^I)\in C_c^2(\Omega_{T_0,R''};\R^N)$ with support contained in $\{r-t< R\}$, we have
\eq{\label{prop3.1con}\tau\int_{\Omega_{T_0,R''}}  [|\partial_t\psi|^2+(\tau|\psi|)^2] t^{-\delta} \ dxdt\lesssim \int_{\Omega_{T_0,R''}} |\wt{\Box}_{\tau,\delta,R'}\psi |^2 t^2\ dxdt. }
Here $\wt{\Box}_{\tau,\delta,R'}$ is defined by\eq{\label{boxfdefn}\wt{\Box}_{\tau,\delta,R'}\psi&:=e^{\tau(f_\delta-R')}\wt{\Box}(e^{-\tau(f_\delta-R')}\psi).}
Here $\wt{\Box}$ is defined by \eqref{lpddefn}.
Note that the constant in \eqref{prop3.1con} is independent of $\tau$, $T_0$, $\delta$ and $\psi$. Also note that $T_0$ can be chosen to be independent of $R'$.
}
\rmk{\rm We  emphasize that here the function $\psi$  is different from the solution $\phi$ in Theorem \ref{mthm}. Thus, in order to apply this proposition, we do  not need to check that $\phi\in C_c^2(\Omega_{T_0,R''};\R^N)$. In fact, in Section \ref{sec4} we will apply Proposition \ref{prop3.1} to the product of $\phi$  and  some cutoff functions.}
\rmk{\rm Let us compare the  estimate \eqref{prop3.1con} with the classical Carleman estimate \eqref{s5carestf1}  for $\wt{\Box}$. Rewrite \eqref{s5carestf1} as follows (by replacing $\psi$ with $\exp(-\tau\exp(\lambda\wt{f}))\psi$):
\eq{\label{clacarest}\tau\norm{\partial\psi}_{L^2(\mcl{N})}^2+\tau^3\norm{\psi}_{L^2(\mcl{N})}^2\sim\tau\norm{\psi}_{H^1_\tau(\mcl{N})}^2&\lesssim \norm{\exp(\tau\exp(\lambda\wt{f}))\wt{\Box}\exp(-\tau\exp(\lambda\wt{f}))\psi}_{L^2(\mcl{N})}^2.}
Note that \eqref{prop3.1con} and \eqref{clacarest} have essentially the same form. However, the integral domain in \eqref{prop3.1con} is a small neighborhood of a point, while that in \eqref{prop3.1con} is not. In other words, \eqref{clacarest} is a local (and thus weaker) estimate while \eqref{prop3.1con} is nonlocal (and thus stronger). To compensate for this, for large time $t\gg1$, we add an extra weight $t^{-\delta}$ to make the left hand side smaller, and an extra weight $t^2$ to make the right hand side larger.
}

\rm\bigskip

Here is a sketch of the proof of Proposition \ref{prop3.1}. In Section \ref{sec3.1}, we decompose the operator $\wt{\Box}_{\tau,\delta,R'}$ into three auxiliary operators $P_1$, $P_2$ and $R^I$. That is, we set
\fm{(\wt{\Box}_{\tau,\delta,R'}\psi)^I=P_1(\psi^I)+P_2(\psi^I)+R^I(\psi),\qquad\forall I}
with
\fm{P_1(\psi^I)&:= \Box\psi^I+\tau^2(\partial^\alpha f\partial_\alpha f )\psi^I,\qquad\text{with symbol }p_1=-m^{\alpha\beta}\xi_\alpha\xi_\beta+\tau^2\partial^\alpha f\partial_\alpha f;}
\fm{P_2(\psi^I)&:=-2\tau\partial^\alpha f\partial_\alpha \psi^I,\qquad\text{with symbol }p_2=-2i\tau\partial^\alpha f \xi_\alpha;}
\fm{R^I(\psi)&:=-\tau(\Box f)\psi^I+V^{I,\alpha}_J\partial_\alpha\psi^J-\tau V_{J}^{I,\alpha}(\partial_\alpha f)\psi^J+W^{I}_J\psi^J.}
In the same subsection, we also define a weighted inner product $\khw{\cdot,\cdot}$ and a corresponding weighted norm $\normw{\cdot}$. The reason why we introduce a weight in the  inner product and the norm will be explained later. With these definitions, we write the right hand side of \eqref{prop3.1con} as 
\fm{\sum_I\normw{(\wt{\Box}_{\tau,\delta,R'}\psi)^I}^2&=\sum_I(\khw{P_1\psi^I,P_1\psi^I}+\khw{P_2\psi^I,P_2\psi^I}+2\khw{P_1\psi^I,P_2\psi^I}+\text{remainders})\\
&=\sum_I(\khw{P_1\psi^I,P_1\psi^I}+\khw{P_2\psi^I,P_2\psi^I}+2\khw{[P_1,P_2]\psi^I,\psi^I}+\text{remainders}).}See \eqref{allest} and \eqref{s3imp} for the accurate formulas. 

Next, in Section \ref{sec3.2}, we handle the terms in \eqref{allest} and \eqref{s3imp} involving the commutator $[P_1,P_2]$. Similar to the proof of other Carleman estimates, a key step in our proof is to estimate the lower bound of the Poisson bracket $\{p_1,p_2\}$; see Lemma \ref{lempois}.   Then, in Section \ref{sec3.3}, we estimate those terms in \eqref{allest} and \eqref{s3imp} involving $R^I$. And finally, in Section \ref{sec3.4}, we conclude our proof of \eqref{prop3.1con}.

\subsection{Setup}  \label{sec3.1}

By the Leibniz's rule, for each $I=1,2,\dots,N$, we can write
\eq{\label{s3tldbox}(\wt{\Box}_{\tau,\delta,R'}\psi)^I=P_1(\psi^I)+P_2(\psi^I)+R^I(\psi),}
where
\eq{P_1(\psi^I)&:= \Box\psi^I+\tau^2(\partial^\alpha f\partial_\alpha f )\psi^I,}
\eq{P_2(\psi^I)&:=-2\tau\partial^\alpha f\partial_\alpha \psi^I,}
\eq{R^I(\psi)&:=-\tau(\Box f)\psi^I+V^{I,\alpha}_J\partial_\alpha\psi^J-\tau V_{J}^{I,\alpha}(\partial_\alpha f)\psi^J+W^{I}_J\psi^J.}
For simplicity, we use the Einstein summation convention. The symbols of $P_1$ and $P_2$ are, respectively,
\eq{\label{p1small}p_1(t,x,\xi,\tau)&=-m^{\alpha\beta}\xi_\alpha\xi_\beta+\tau^2\partial^\alpha f\partial_\alpha f,}
and
\eq{\label{p2small}p_2(t,x,\xi,\tau)&=-2i\tau\partial^\alpha f \xi_\alpha.}
Here we assume $(t,x)\in\R^{1+3}$ and $\xi=(\xi_\alpha)_{\alpha=0}^3\in\R^{1+3}$. We notice that $p_1$ and $p_2$ are both homogeneous polynomials of $(\tau,\xi)\in\R^{1+1+3}$ of order $2$, that $p_1$ is real and $p_2$ is purely imaginary, and that the symbol of the $R^I$'s are polynomials of $(\tau,\xi)\in\R^{1+1+3}$ of order $1$. This explains how we make the decomposition \eqref{s3tldbox}.

For two $\R$-valued functions $F=F(t,x)$ and $G=G(t,x)$, we write \eq{\khw{F,G}:=\int_{\Omega_{T_0,R''}}FGt^2\ dxdt}
and
\eq{\normw{F}:=\sqrt{\khw{F,F}}.}
The reason why we add a weight will be explained later; see Remark \ref{whyweighted}.

We can  compute the conjugate under this inner product.
\lem{\label{lemconj}Suppose that $F,G\in C^1_c(\Omega_{T_0,R''})$. Then, we have 
\fm{\khw{P_j(F),G}=\khw{F,P_j^*(G)},\qquad j=1,2,}
where
\eq{\label{p1star}P_1^*&:=P_1-2t^{-2}-4t^{-1}\partial_t ,}and
\eq{\label{p2star}P_2^*&:=-P_2+2\tau \Box f-4\tau t^{-1}\partial_tf.}}
\begin{proof}Fix $F,G$ as in the statement. By integration by parts, we have
\fm{\int_{\Omega_{T_0,R''}}(\Box F)Gt^2\ dxdt&=\int_{\Omega_{T_0,R''}}F\Box (Gt^2)\ dxdt=\int_{\Omega_{T_0,R''}}F (t^2\Box G-4t\partial_tG-2G)\ dxdt,}
\fm{\int_{\Omega_{T_0,R''}}(-2\tau \partial^\alpha f\partial_\alpha F)Gt^2\ dxdt&=\int_{\Omega_{T_0,R''}}2\tau F\partial_\alpha((\partial^\alpha f) Gt^2)\ dxdt\\
&=\int_{\Omega_{T_0,R''}}2\tau F((\Box f) Gt^2+t^2\partial^\alpha f\partial_\alpha G-2Gt\partial_t f)\ dxdt.}
Then \eqref{p1star} and \eqref{p2star} follow.
\end{proof}
\rmk{\label{lemconjrmk}\rm Suppose that $F,G\in C^1_c(\Omega_{T_0,R''})$.  For general $C^1$ vector field $X=X^\alpha \partial_\alpha$, we have
\fm{\khw{X(F),G}=\khw{F,-t^{-2}\partial_\alpha(t^2X^\alpha G)}.}
This identity follows easily from integration by parts.}\rm\bigskip

As a result, the right side of \eqref{prop3.1con} is equal to \eq{\label{allest}\sum_{I=1}^N\normw{(\wt{\Box}_{\tau,\delta,R'}\psi)^I}^2&=\sum_I\left[\normw{P_1(\psi^I)}^2+\normw{P_2(\psi^I)+R^I(\psi)}^2\right.\\&\hspace{5em}\left.+\khw{(P_2^*P_1+P_1^*P_2)(\psi^I),\psi^I}+2\khw{P_1(\psi^I),R^I(\psi)}\right].}
By \eqref{p1star} and \eqref{p2star}, we have
\fm{(P_2^*P_1+P_1^*P_2)(\psi^I)&=[P_1,P_2]\psi^I+(2\tau \Box f-4\tau t^{-1}\partial_tf)P_1(\psi^I)+(-2t^{-2}-4t^{-1}\partial_t )P_2(\psi^I).}
By applying Remark \ref{lemconjrmk} to $X=\partial_t$, we have
\eq{\label{s3imp}&\khw{(P_2^*P_1+P_1^*P_2)(\psi^I),\psi^I}+2\khw{P_1(\psi^I),R^I(\psi)}\\&=\khw{[P_1,P_2]\psi^I,\psi^I}+2\khw{P_1(\psi^I),(\tau \Box f-2\tau t^{-1}\partial_tf)\psi^I+R^I(\psi)}\\
&\quad-\khw{P_2(\psi^I),2t^{-2}\psi^I}+\khw{P_2(\psi^I),4t^{-2}\partial_t(t\psi^I)}.}
A key step in the proof is to obtain a lower bound for \eqref{s3imp}.

We end this setup with two remarks involving the last term, $\khw{P_2(\psi^I),4t^{-2}\partial_t(t\psi^I)}$, in \eqref{s3imp}. First, this term appears because of the  weight in $\khw{\cdot,\cdot}$. If we replace $\khw{\cdot,\cdot}$ with an unweighted inner product $\int FG\ dxdt$, then we will obtain a similar identity but without a term like $\khw{P_2(\psi^I),4t^{-2}\partial_t(t\psi^I)}$. Second, it turns out that the term $\khw{P_2(\psi^I),4t^{-2}\partial_t(t\psi^I)}$ is crucial in the proof of the Carleman estimate \eqref{prop3.1con}. In fact, it will cancels with another term from the next proposition (Proposition \ref{prop3.1.1}) which we cannot estimate directly. See Remark \ref{whyweighted} below.

\subsection{Estimates for $\khw{[P_1,P_2]\psi^I,\psi^I}$} \label{sec3.2}
Our goal now is to estimate $\khw{[P_1,P_2]\psi^I,\psi^I}$. In fact we have the following proposition. \prop{\label{prop3.1.1}Fix $T_0\gg_{R,R''}1$. For $\tau\gg1$, we have
\eq{\label{prop3.1.1f1}&\khw{[P_1,P_2]\psi^I,\psi^I}\\&\geq -\tau^{-1}\normw{P_2\psi^I}^2+ \int_{\Omega_{T_0,R''}}  \frac{1}{4}\tau t^{-\delta}[(\partial_t\psi^I)^2+\tau^2  (\psi^I)^2]\ dxdt\\&\quad-4\tau\khw{P_1\psi^I,r^{-1}\psi^I}+4\khw{ P_2\psi^I,r^{-1}(-1+t^{-\delta-1})\partial_t\psi^I}.}}
\rmk{\label{whyweighted}\rm The last term, $4\khw{ P_2\psi^I,r^{-1}(-1+t^{-\delta-1})\partial_t\psi^I}$, is the reason why we need to use the weighted inner product $\khw{\cdot,\cdot}$ and weighted norm $\normw{\cdot}$. In fact, this term will always appear, no matter whether we use the weighted inner product or the unweighted one. We also notice that this term cannot be controlled by applying the Cauchy-Schwarz directly, because we will get a norm
\fm{\tau\normw{r^{-1}(-1+t^{-\delta-1})\partial_t\psi^I}\sim\tau\int (\partial_t\psi^I)^2\ dxdt}
which cannot be absorbed by the positive integral in \eqref{prop3.1.1f1}. Similarly for the unweighted case. The only way to control $4\khw{ P_2\psi^I,r^{-1}(-1+t^{-\delta-1})\partial_t\psi^I}$ is to cancel it with another term, and there is a cancellation in the weighted case (also see the remark after \eqref{s3imp})
\fm{4\khw{ P_2\psi^I,r^{-1}(-1+t^{-\delta-1})\partial_t\psi^I}+\khw{P_2(\psi^I),4t^{-2}\partial_t(t\psi^I)}&=\khw{P_2\psi^I,\text{lower order terms}},}
where the norm of the lower order terms can be absorbed by the positive integral in \eqref{prop3.1.1f1}, so we can use the Cauchy-Schwarz inequality to estimate it. This explains why we need the weight in the definition of $\khw{\cdot,\cdot}$.
}
\rm\bigskip

We now  prove Proposition \ref{prop3.1.1}.  The proof relies on the lower bound of the symbol of $[P_1,P_2]$ and integration by parts. To estimate $\khw{[P_1,P_2]\psi^I,\psi^I}$, we need to first compute the Poisson bracket $\{p_1,p_2\}$ where $p_j$ is the symbol of $P_j$ for $j=1,2$. 

\lem{\label{lempois} Fix $T_0\gg_{R,R''} 1$. For each $(t,x)\in\Omega_{T_0,R''}$ such that $r-t<R$ and for each $\xi\in\R^{1+3}$, we have
\eq{\label{poiest}\frac{1}{i\tau}\{p_1,p_2\}&\geq  -r^{-1}(\frac{p_2}{i\tau})^2-4 r^{-1}p_1+  t^{-\delta-2}(\xi_0^2+\tau^2)+4r^{-1}\cdot\frac{p_2}{i\tau}\cdot (-1+t^{-\delta-1})\xi_0.}
}
\begin{proof} We have
\fm{\frac{1}{i\tau}\{p_1,p_2\}&=-2m^{\alpha\beta}\xi_\beta \partial_\alpha(-2 \partial^\nu f\xi_\nu)-(-2 \partial^\alpha f)\partial_\alpha (\tau^2\partial^\beta f\partial_\beta f)\\
&=4m^{\alpha\beta}\xi_\beta \xi_\nu \partial_\alpha\partial^\nu f+2\tau^2 \partial^\alpha f \partial_\beta f \partial_\alpha  \partial^\beta f +2\tau^2 \partial^\alpha f \partial^\beta f \partial_\alpha\partial_\beta f\\
&=4\xi_\beta \xi_\nu \partial^\beta\partial^\nu f+4\tau^2 \partial_\alpha f \partial_\beta f \partial^\alpha\partial^\beta f.}
Note that $(\partial_\alpha f)_{\alpha=0}^3=(-1+t^{-\delta-1},\omega_1,\omega_2,\omega_3)$ and $(\partial^\alpha f)_{\alpha=0}^3=(1-t^{-\delta-1},\omega_1,\omega_2,\omega_3)$. Thus, $\partial^\alpha f\partial_\alpha f=1-(1-t^{-\delta-1})^2>0$. Moreover,
\fm{\partial^i\partial^j f&=r^{-1}(\delta_{ij}-\omega_i\omega_j),\qquad i,j=1,2,3;\\\partial^0\partial^0 f&=-(\delta+1)t^{-\delta-2};\\
\partial^0\partial^i f&=\partial^i\partial^0 f=0,\qquad i=1,2,3.}
Thus, $\xi_\alpha\xi_\beta\partial^\alpha\partial^\beta f=-(\delta+1)t^{-\delta-2}\xi_0^2+4r^{-1}(\sum\xi_i^2-(\sum \omega_i\xi_i)^2)$. It then follows that
\fm{\frac{1}{i\tau}\{p_1,p_2\}&=-4(\delta+1)t^{-\delta-2}\xi_0^2 +4r^{-1}(\sum_{i=1}^3\xi_i^2-(\sum_{i=1}^3\omega_i\xi_i)^2)-4\tau^2(\delta+1)t^{-\delta-2}(-1+t^{-\delta-1})^2 .}
Moreover, by \eqref{p1small} and \eqref{p2small},  we have
\fm{p_1&=\xi_0^2-\sum_{i=1}^3\xi_i^2+\tau^2t^{-\delta-1}(2-t^{-\delta-1}),\\
p_2&=-2i\tau ((1-t^{-\delta-1})\xi_0+\sum_{i=1}^3\omega_i\xi_i).}
Using $p_1$ and $p_2$, we can express $\sum_{i=1}^3\xi_i^2-(\sum_{i=1}^3\omega_i\xi_i)^2$ in terms of $\xi_0$. That is,
\fm{&\sum_{i=1}^3\xi_i^2-(\sum_{i=1}^3\omega_i\xi_i)^2\\&=\xi_0^2-p_1+\tau^2t^{-\delta-1}(2-t^{-\delta-1})-(\frac{p_2}{-2i\tau}-(1-t^{-\delta-1})\xi_0)^2\\
&=t^{-\delta-1}(2-t^{-\delta-1})\xi_0^2-p_1+\tau^2t^{-\delta-1}(2-t^{-\delta-1})-(\frac{p_2}{2i\tau})^2-\frac{p_2}{i\tau}(1-t^{-\delta-1})\xi_0,}
and
\eq{\label{lempoisf1}\frac{1}{i\tau}\{p_1,p_2\}&=-4r^{-1}p_1-r^{-1}(\frac{p_2}{i\tau})^2-4r^{-1}\frac{p_2}{i\tau}(1-t^{-\delta-1})\xi_0\\&\quad+4r^{-1}t^{-\delta-2}(-(\delta+1) r+t(2-t^{-\delta-1}))\xi_0^2\\&\quad+4\tau^2r^{-1}t^{-\delta-2}(t(2-t^{-\delta-1})-(\delta+1)r(-1+t^{-\delta-1})^2) .}
Since $0<\delta<1/2$, $R''<r-t<R$ and $t\geq T_0$, we have
\fm{-(\delta+1) r+t(2-t^{-\delta-1})\geq 2t-\frac{3}{2}r-t^{-\delta}\geq 2t-\frac{3}{2}(t+R)-T_0^{-\delta}\geq \frac{1}{2}t-\frac{3}{2}|R|-T_0^{-\delta}.}
If $T_0\gg_{R}1$ (say $T_0\geq 6+9|R|$), we have $t\geq 9|R|+6$ and thus \fm{-(\delta+1) r+t(2-t^{-\delta-1})\geq  \frac{1}{2}t-\frac{3}{2}|R|-1\geq \frac{1}{3}t.}
Moreover, we have
\fm{&t(2-t^{-\delta-1})-(\delta+1)r(-1+t^{-\delta-1})^2\geq t(2-t^{-\delta-1})-(\delta+1)r}
which is also no less than $t/3$ as proved above. Moreover, for $T_0\gg_{R,R''}1$, since $R''<r-t<R$, we have $|r-t|<t/3$ and thus $r/4<t/3$. By combining these  estimates with \eqref{lempoisf1},  we obtain \eqref{poiest}.
\end{proof}
\rmk{\rm If we set $f=r-t$ and do the same computations, we will get
\fm{\frac{1}{i\tau}\{p_1,p_2\}&=4r^{-1}(\sum_{i=1}^3\xi_i^2-(\sum_{i=1}^3\omega_i\xi_i)^2)=-4r^{-1}p_1-r^{-1}(\frac{p_2}{i\tau})^2-4r^{-1}\frac{p_2}{i\tau}\xi_0.}
Comparing this with \eqref{poiest}, we notice that we lose a positive term $t^{-\delta-2}(\xi_0^2+\tau^2)$. Here we obtain such a positive term in \eqref{poiest}  because of the term $-\delta^{-1} t^{-\delta}$ in $f$ which makes $\partial^0f=-\partial_tf<1$ and $\partial^0\partial^0f=\partial_t^2f<0$. This is the key in the proof of the Carleman estimate.}

\rm\bigskip

Set
\fm{b(t,x,\xi,\tau):=\frac{1}{i\tau}\{p_1,p_2\}=4\xi_\beta \xi_\nu \partial^\beta\partial^\nu f+4\tau^2 \partial_\alpha f \partial_\beta f \partial^\alpha\partial^\beta f.}
It is related to the following differential operator
\fm{b(t,x,D,\tau)\Phi:=-4 \partial^\beta\partial^\nu f \partial_\beta\partial_\nu \Phi+4\tau^2 \partial_\alpha f \partial_\beta f \partial^\alpha\partial^\beta f\cdot\Phi.}
So far, we have proved a lower bound for the symbol $b$. The next lemma  allows us to use this lower bound to control $\khw{ [P_1,P_2]\psi^I,\psi^I}$.

\lem{\label{bctrl} We have
\eq{\label{bctrlf1}\khw{ [P_1,P_2]\psi^I,\psi^I}&=\int_{\Omega_{T_0,R''}}\tau b(t,x,\partial\psi^I,\tau\psi^I)t^2\ dxdt+2\tau\khw{\psi^I,(\partial^\nu\Box f)\partial_\nu\psi^I}\\&\quad-8\tau(\delta+1)\khw{t^{-\delta-3}\psi^I,\partial_t\psi^I}.}}
\begin{proof}We notice that \fm{\ [P_1,P_2]&=\tau b(t,x,D,\tau)+\Box(-2\tau\partial^\alpha f)\partial_\alpha=\tau b(t,x,D,\tau)-2\tau(\partial^\alpha\Box f)\partial_\alpha }
and therefore
\fm{\khw{[P_1,P_2]\psi^I,\psi^I}&=\khw{\tau b(t,x,D,\tau)\psi^I,\psi^I}-2\tau\khw{(\partial^\alpha\Box f)\partial_\alpha\psi^I,\psi^I}.}
In addition, since
\fm{&\int_{\Omega_{T_0,R''}} -4 \partial^\beta\partial^\nu f \partial_\beta\partial_\nu\psi^I\cdot\psi^It^2\ dxdt=\int_{\Omega_{T_0,R''}} 4\partial_\nu\psi^I\cdot \partial_\beta(\partial^\beta\partial^\nu f \cdot\psi^It^2)\ dxdt\\
&=\int_{\Omega_{T_0,R''}} 4\partial_\nu\psi^I\cdot \partial_\beta(t^2\partial^\beta\partial^\nu f) \cdot\psi^I+ 4\partial_\nu\psi^I\cdot \partial_\beta \psi^I\cdot t^2\partial^\beta\partial^\nu f \ dxdt,}
we have
\fm{\khw{\tau b(t,x,D,\tau)\psi^I,\psi^I}&=\int_{\Omega_{T_0,R''}} 4\tau\partial_\nu\psi^I\cdot \partial_\beta(t^2\partial^\beta\partial^\nu f) \cdot\psi^I+\tau b(t,x,\partial\psi^I,\tau\psi^I) t^2 \ dxdt.}
We also notice that
\fm{4\tau\partial_\nu\psi^I\cdot \partial_\beta(t^2\partial^\beta\partial^\nu f)\psi^I&=4\tau\psi^I\partial_\nu\psi^I\cdot t^2\partial^\nu\Box f+4\tau\psi^I\partial_t\psi^I\cdot 2t\partial^0\partial^0 f\\
&=4\tau\psi^I\partial_\nu\psi^I\cdot t^2\partial^\nu\Box f-8(\delta+1)\tau t^{-\delta-1}\psi^I\partial_t\psi^I,}
so
\fm{&\int_{\Omega_{T_0,R''}}4\tau \partial_\nu\psi^I\cdot \partial_\beta(t^2\partial^\beta\partial^\nu f) \cdot\psi^I\ dxdt\\&=4\tau\khw{\psi^I,(\partial^\nu\Box f)\partial_\nu\psi^I}-8(\delta+1)\tau\int_{\Omega_{T_0,R''}}t^{-\delta-1}\psi^I\partial_t\psi^I\ dxdt.}
Then \eqref{bctrlf1} follows.
\end{proof}
\rm\bigskip

Let us estimate the first term on the right side of \eqref{bctrlf1} in Lemma \ref{bctrl}.
\lem{\label{lemcomctrl} Fix $T_0\gg_{R,R''}1$. For $\tau\gg1$, we have
\eq{\label{lemcomctrlf1} &\int_{\Omega_{T_0,R''}}\tau b(t,x,\partial\psi^I,\psi^I)t^2\ dxdt\\&\geq -\tau^{-1}\normw{P_2\psi^I}^2+ \int_{\Omega_{T_0,R''}}  \tau t^{-\delta}[(\partial_t\psi^I)^2+\frac{1}{2}\tau^2  (\psi^I)^2]\ dxdt\\&\quad-4\tau\khw{P_1\psi^I,r^{-1}\psi^I}+4\khw{ P_2\psi^I,r^{-1}(-1+t^{-\delta-1})\partial_t\psi^I}+\int_{\Omega_{T_0,R''}} 4\tau r^{-2}t^2\psi^I \partial_r\psi^I \ dxdt.}}
\begin{proof}
Note that $(p_2/(i\tau))(t,x,\partial\psi^I,\tau\psi^I)=\tau^{-1}P_2\psi^I$. So by Lemma \ref{lempois}, we have
\eq{\label{lemcomctrlf2}&\int_{\Omega_{T_0,R''}}\tau b(t,x,\partial\psi^I,\psi^I)t^2\ dxdt\\&\geq \int_{\Omega_{T_0,R''}}   \tau t^{-\delta}[(\partial_t\psi^I)^2+\tau^2  (\psi^I)^2]-r^{-1}t^2\tau^{-1}(P_2\psi^I)^2\\&\qquad\qquad-4\tau r^{-1}t^2p_1(t,x,\partial\psi^I,\tau\psi^I)+4r^{-1}t^2 P_2\psi^I(-1+t^{-\delta-1})\partial_t\psi^I\ dxdt\\
&\geq -\tau^{-1}\normw{P_2\psi^I}^2+ \int_{\Omega_{T_0,R''}}  \tau  t^{-\delta}[(\partial_t\psi^I)^2+\tau^2  (\psi^I)^2]\ dxdt\\&\quad+\int_{\Omega_{T_0,R''}} -4\tau r^{-1}t^2p_1(t,x,\partial\psi^I,\tau\psi^I)\ dxdt+4\khw{ P_2\psi^I,r^{-1}(-1+t^{-\delta-1})\partial_t\psi^I} .}
Here we recall that $\psi^I$ is nonzero only if $R''<r-t<R$ and $t\geq T_0$. Since $T_0\gg_{R,R''}1$, we can always assume $r>t+R''\geq T_0-|R''|\geq 4$.  Now we need to compute \fm{\int_{\Omega_{T_0,R''}} -4\tau r^{-1}t^2p_1(t,x,\partial\psi^I,\tau\psi^I)\ dxdt.} It is equal to
\fm{&\int_{\Omega_{T_0,R''}} -4\tau r^{-1}t^2(-m^{\alpha\beta}\partial_\alpha\psi^I\partial_\beta\psi^I+\tau^2t^{-\delta-1}(2-t^{-\delta-1})(\psi^I)^2)\ dxdt\\
&=\int_{\Omega_{T_0,R''}} -4\tau\psi^I\partial_\alpha( r^{-1}t^2 m^{\alpha\beta}\partial_\beta\psi^I) \ dxdt+\int_{\Omega_{T_0,R''}} -4\tau^3 r^{-1} t^{1-\delta}(2-t^{-\delta-1})(\psi^I)^2\ dxdt\\
&=-4\tau\khw{P_1\psi^I,r^{-1}\psi^I}+\int_{\Omega_{T_0,R''}} -4\tau\psi^I\partial^\beta( r^{-1}t^2) \partial_\beta\psi^I \ dxdt\\
&=-4\tau\khw{P_1\psi^I,r^{-1}\psi^I}+\int_{\Omega_{T_0,R''}} 8\tau r^{-1}t \psi^I  \partial_t\psi^I+4\tau r^{-2}t^2\psi^I \partial_r\psi^I \ dxdt.}
By integration by parts again, we have
\fm{\int_{\Omega_{T_0,R''}} 8\tau r^{-1}t \psi^I  \partial_t\psi^I \ dxdt&=-\int_{\Omega_{T_0,R''}} 8\tau r^{-1}\partial_t(t \psi^I ) \psi^I \ dxdt\\
&=-\int_{\Omega_{T_0,R''}} 8\tau r^{-1}(t \partial_t\psi^I+\psi^I ) \psi^I \ dxdt}
and thus
\fm{\int_{\Omega_{T_0,R''}} 8\tau r^{-1}t \psi^I  \partial_t\psi^I \ dxdt
&=-\int_{\Omega_{T_0,R''}} 4\tau r^{-1}(\psi^I )^2  \ dxdt\\&\geq -4\tau^{-1}T_0^{\delta}(T_0-|R|-|R''|)^{-1}\int_{\Omega_{T_0,R''}} \tau^3  t^{-\delta}(\psi^I )^2  \ dxdt.}
By choosing  $T_0\gg_{R,R''} 1$, we have $T_0^{\delta}(T_0-|R|-|R''|)^{-1}\leq 1$. As long as $\tau>8$, we can finish the proof since in \eqref{lemcomctrlf2}  there is a term \fm{\int_{\Omega_{T_0,R''}} \tau t^{-\delta}[(\partial_t\psi^I)^2+\tau^2  (\psi^I)^2] \ dxdt.}
\end{proof}\rm\bigskip

Lemma \ref{lemcomctrl} gives us a lower bound for the first term in the expansion \eqref{bctrlf1} of $\khw{ [P_1,P_2]\psi^I,\psi^I}$. To finish the proof of Proposition \ref{prop3.1.1}, we need to estimate the remaining two terms on the right side of \eqref{bctrlf1}.
Since 
\fm{\Box f&=(\delta+1)t^{-\delta-2}+2r^{-1},}
we have
\fm{&2\tau\khw{\psi^I,(\partial^\nu\Box f)\partial_\nu\psi^I}=2(\delta+1)(\delta+2)\tau\khw{\psi^I,t^{-\delta-3} \partial_t\psi^I}-4\tau\khw{\psi^I,r^{-2} \partial_r\psi^I}.}
Recall that in  \eqref{lemcomctrlf1}, we have 
\fm{\int_{\Omega_{T_0,R''}}\tau b(t,x,\partial\psi^I,\psi^I)t^2\ dxdt&\geq\text{other terms}+\int_{\Omega_{T_0,R''}} 4\tau r^{-2}t^2\psi^I \partial_r\psi^I \ dxdt\\
&=\text{other terms}+4\tau\khw{\psi^I,r^{-2} \partial_r\psi^I},}
so  $4\tau\khw{\psi^I,r^{-2} \partial_r\psi^I}$ cancels with the last term in \eqref{lemcomctrlf1}. In addition, by the Cauchy-Schwarz inequality, we have
\fm{&2(\delta+1)(\delta+2)\tau\khw{\psi^I,t^{-\delta-3} \partial_t\psi^I}-8\tau(\delta+1)\khw{t^{-\delta-3}\psi^I,\partial_t\psi^I}\\&\geq -CT_0^{-1}\tau\normw{t^{-1-\delta/2}\psi^I}\normw{t^{-1-\delta/2}\partial_t\psi^I}\\
&\geq -C\tau^{-1}\cdot\tau(\tau^2\normw{t^{-1-\delta/2}\psi^I}^2+\normw{t^{-1-\delta/2}\partial_t\psi^I}^2).}
Note that $\tau(\tau^2\normw{t^{-1-\delta/2}\psi^I}^2+\normw{t^{-1-\delta/2}\partial_t\psi^I}^2)$ appears on the right side of \eqref{lemcomctrlf1}. Thus, we finish the proof by choosing $\tau\gg 1$ such that $C\tau^{-1}<1/8$.

\subsection{Estimates for $\normw{R^I(\psi)+2\tau r^{-1}\psi^I}$ and $\normw{P_2(\psi^I)+R^I(\psi)}$}\label{sec3.3} Recall that 
\fm{R^I(\psi)&=-\tau(\Box f)\psi^I+V^{I,\alpha}_J\partial_\alpha\psi^J-\tau V_{J}^{I,\alpha}(\partial_\alpha f)\psi^J+W^{I}_J\psi^J\\
&=-2\tau r^{-1}\psi^I-\tau(\delta+1)t^{-\delta-2}\psi^I+V^{I,\alpha}_J\partial_\alpha\psi^J-\tau V_{J}^{I,\alpha}(\partial_\alpha f)\psi^J+W^{I}_J\psi^J.}
The main proposition in this subsection is the following.
\prop{\label{prop3.7} For $T_0\gg_{R,R'',\gamma,C_1}1$ and for $\tau\gg1$, we have
\eq{\label{prop3.7f1}&\normw{R^I(\psi)+2\tau r^{-1}\psi^I}\\&\leq T_0^{-\gamma/4}\sum_J(\normw{P_1\psi^J}+\normw{P_2\psi^J}+\tau\normw{t^{-\delta/2-1}\psi^J}+\normw{t^{-\delta/2-1}\partial_t\psi^J})+2\tau^{-1/2}\sum_J\normw{P_2\psi^J}}and
\eq{\label{prop3.7f2}&\normw{P_2(\psi^I)+R^I(\psi)}^2\\&\geq \frac{1}{100}\normw{P_2\psi^I}^2-\frac{1}{1000N}\sum_{J}(\normw{P_1\psi^J}^2+\normw{P_2\psi^J}^2+\tau\normw{t^{-\delta/2-1}\psi^J}^2+\normw{t^{-\delta/2-1}\partial_t\psi^J}^2).}
}\rm\bigskip

Here is a sketch of the proof. Using the pointwise bounds \eqref{prop3.1a1} and \eqref{prop3.1a2} for $V^{*,*}_*,W^*_*$, we first show that
\eq{\label{prop3skpf}\normw{R^I(\psi)+2\tau r^{-1}\psi^I}&\leq C\tau \normw{t^{-\delta-2}\psi^I}+C_1(2\tau+1)\normw{(1+t+r)^{-1-\gamma}\psi}+\normw{V^{I,\alpha}_J\partial_\alpha\psi^J}.}
The first two terms on the right hand side can be controlled by $T_0^{-\gamma/2}\normw{P_2\psi}$ if we  apply a Hardy-type inequality \eqref{lem3.3.1f1}. To estimate $V^{I,\alpha}_J\partial_\alpha\psi^J$, we write \fm{V_J^{I,\alpha}\partial_\alpha\psi^J&=\frac{1}{2}V^{I,\alpha}_J\wt{\omega}_\alpha(\partial_t+\partial_r)\psi^J+\frac{1}{2}V^{I,\alpha}_J\widehat{\omega}_\alpha(-\partial_t+\partial_r)\psi^J+X_\omega\psi^J.}
Here $X_\omega$ is a linear combination of $(\partial_i-\omega_i\partial_r)_{i=1,2,3}$. By writing $P_2=-2\tau(\partial_t+\partial_r)+\text{lower order terms}$, we can control the first term by $P_2\psi$ and some remainders. Because of  \eqref{prop3.1a2}, the coefficient of $(-\partial_t+\partial_r)\psi^J$ is $O((1+t+r)^{-1-\gamma})$ and we can control the second term directly. The last term can be controlled by applying Lemma \ref{lem3.3.2}. Combining all these results, we obtain \eqref{prop3.7f1}. The estimate \eqref{prop3.7f2} now follows from  \eqref{prop3.7f1} and \eqref{lem3.3.1f1}.

We start with the following lemma.
\lem{\label{lem3.3.1} For $T_0\gg_{R,R''}1$, we have\eq{\label{lem3.3.1f1}\normw{r^{-1}\psi^I}\leq\frac{2}{5\tau}\normw{P_2\psi^I}.}}
\begin{proof}
By \eqref{p2star}, we have
\fm{P_2+P_2^*&=2\tau (2r^{-1}+(\delta+1)t^{-\delta-2})-4\tau t^{-1}(-1+t^{-\delta-1})\\
&=8\tau r^{-1}+2\tau (\delta-1)t^{-\delta-2}+4\tau (rt)^{-1}(r-t).}
Then,
\fm{&\int_{\Omega_{T_0,R''}} r^{-2}t^{2}(\psi^I)^2\ dxdt\\
&=\int_{\Omega_{T_0,R''}} r^{-1}\psi^I\cdot t^{2}\cdot\frac{P_2\psi^I+P_2^*\psi^I+[2\tau(1-\delta)t^{-\delta-2}-4\tau (rt)^{-1}(r-t)]\psi^I}{8\tau}\ dxdt\\
&=\frac{1}{8\tau}(\khw{r^{-1}\psi^I,P_2\psi^I}+\khw{P_2(r^{-1}\psi^I),\psi^I})\\&\quad+\int_{\Omega_{T_0,R''}}\frac{1}{4}((1-\delta)r^{-1}t^{-\delta}-2 r^{-2}t(r-t))(\psi^I)^2\ dxdt\\
&= \frac{1}{4\tau}\khw{r^{-1}\psi^I,P_2\psi^I}+\frac{1}{4}\khw{r^{-2}\psi^I,\psi^I}\\&\quad+\int_{\Omega_{T_0,R''}}(\frac{1}{4}(1-\delta)rt^{-\delta-2}-\frac{1}{2} t^{-1}(r-t))\cdot r^{-2}t^2(\psi^I)^2\ dxdt\\
&\leq \frac{1}{4\tau}\normw{r^{-1}\psi^I}\normw{P_2\psi^I}+\frac{1}{4}\normw{r^{-1}\psi^I}^2\\&\quad+\int_{\Omega_{T_0,R''}}(\frac{1}{4}(T_0+|R|+|R''|)T_0^{-\delta-2}+\frac{1}{2} T_0^{-1}(|R|+|R''|))\cdot r^{-2}t^2(\psi^I)^2\ dxdt.}
Note that $r\geq t+R''$ and $t\geq T_0$ whenever $\psi^I\neq 0$. By choosing $T_0\gg_{R,R''}1$, we have
\fm{\frac{1}{4}(T_0+|R|+|R''|)T_0^{-\delta-2}+\frac{1}{2} T_0^{-1}(|R|+|R''|)\leq \frac{1}{8}.}Then,
\fm{\frac{3}{4}\normw{r^{-1}\psi^I}^2\leq\frac{1}{4\tau}\normw{r^{-1}\psi^I}\normw{P_2\psi^I}+\frac{1}{8} \normw{r^{-1}\psi^I}^2. }
Then \eqref{lem3.3.1f1} follows.
\end{proof}
\rm \bigskip

Using the previous lemma, we  estimate $-\tau(\delta+1)t^{-\delta-2}\psi^I-\tau V_{J}^{I,\alpha}(\partial_\alpha f)\psi^J+W^{I}_J\psi^J$ which appears in the expansion of $R^I(\psi)+2\tau r^{-1}\psi^I$. By \eqref{prop3.1a1} and \eqref{prop3.1a2}, for fixed $I$ and $J$, in $\Omega_{T_0,R''}$ we have $|W^I_J|\leq C_1(1+t+r)^{-1-\gamma}$ and
\fm{|V^{I,\alpha}_J\partial_\alpha f|&=|V^{I,\alpha}_J\partial_\alpha (r-t)+V^{I,0}_J\cdot t^{-\delta-1}|\\&\leq C_1(1+t+r)^{-1-\gamma}+C_1(1+t+r)^{-2-\delta}\leq 2C_1(1+t+r)^{-1-\gamma}.}  As a result, we have
\fm{&\normw{-\tau(\delta+1)t^{-\delta-2}\psi^I-\tau V_{J}^{I,\alpha}(\partial_\alpha f)\psi^J+W^{I}_J\psi^J}\\
&\leq \frac{3}{2}\tau \normw{t^{-\delta-2}\psi^I}+C_1(2\tau+1)\sum_J\normw{(1+t+r)^{-1-\gamma}\psi^J}.}
In addition, by Lemma \ref{lem3.3.1} we have
\fm{\frac{3}{2}\tau \normw{t^{-\delta-2}\psi^I}&\leq \frac{3}{2}\tau T_0^{-2-\delta}(T_0+|R|+|R''|)\normw{r^{-1}\psi^I}\leq CT_0^{-2-\delta}(T_0+|R|+|R''|)\normw{P_2\psi^I},}
\fm{C_1(2\tau+1)\tau\normw{(1+t+r)^{-1-\gamma}\psi^J}&\leq C_1(2\tau+1)\tau T_0^{-\gamma} \normw{r^{-1}\psi^J}\leq CC_1T_0^{-\gamma}\normw{P_2\psi^J}.}
By choosing $T_0\gg_{R,R'',\gamma,C_1}1$, we conclude that
\eq{\label{riestf1}&\normw{-\tau(\delta+1)t^{-\delta-2}\psi^I-\tau V_{J}^{I,\alpha}(\partial_\alpha f)\psi^J+W^{I}_J\psi^J}\leq T_0^{-\gamma/2}\sum_J\normw{P_2\psi^J}.}

It remains to estimate $V^{I,\alpha}_J\partial_\alpha\psi^J$. For each fixed $I,J$, we have
\eq{\label{decvjp}V_J^{I,\alpha}\partial_\alpha\psi^J&=\frac{1}{2}V^{I,\alpha}_J\wt{\omega}_\alpha(\partial_t+\partial_r)\psi^J+\frac{1}{2}V^{I,\alpha}_J\widehat{\omega}_\alpha(-\partial_t+\partial_r)\psi^J+X_\omega\psi^J.}
Here $\wt{\omega}=(1,\omega)$, $\widehat{\omega}=(-1,\omega)$, and $X_\omega$ is an angular derivative. Here we need the following lemma which gives an estimate for the $X_\omega$.
\lem{\label{lem3.3.2} Let $\slashed\partial=(\partial_j-\omega_j\partial_r)_{j=1,2,3}$ denote the angular derivatives.  Then, for $\tau\gg1$ we have
\eq{\label{lem3.3.2f1}&\normw{r^{-1}\slashed\partial\psi^I}^2\\&\lesssim T_0^{-2}\normw{P_1\psi^I}^2+\tau^{-1} \normw{P_2\psi^I}^2+\tau^2T_0^{-1}\normw{t^{-\delta/2-1}\psi^I}^2+T_0^{-1}\normw{t^{-\delta/2-1}\partial_t\psi^I}^2.}}
\begin{proof}Since $t\sim r$ in the support of $\psi$, we can estimate $\normw{t^{-1}\slashed\partial\psi^I}$ instead. We have\fm{|\slashed\partial\psi^I|^2&=\sum_{j=1}^3(\partial_j\psi^I)^2-(\partial_r\psi^I)^2=m^{\alpha\beta}\partial_\alpha\psi^I\partial_\beta\psi^I+(\partial_t\psi^I)^2-(-\frac{1}{2\tau}P_2\psi^I-(1-t^{-\delta-1})\partial_t\psi^I)^2\\
&=-p_1(t,x,\partial\psi^I,\tau\psi^I)-(4\tau^2)^{-1}(P_2\psi^I)^2-\tau^{-1}(1-t^{-\delta-1})\partial_t\psi^IP_2\psi^I\\&\quad+\tau^2 t^{-\delta-1}(2-t^{-\delta-1})(\psi^I)^2+t^{-\delta-1}(2-t^{-\delta-1})(\partial_t\psi^I)^2.}
Integrate both sides over $\Omega_{T_0,R''}$. Note that
\fm{\int_{\Omega_{T_0,R''}}-p_1(t,x,\partial\psi^I,\tau\psi^I) \ dxdt&=\int_{\Omega_{T_0,R''}}m^{\alpha\beta}\partial_\alpha\psi^I\partial_\beta\psi^I-\tau^2t^{-\delta-1}(2-t^{-\delta-1})(\psi^I)^2 \ dxdt\\
&=\int_{\Omega_{T_0,R''}}-m^{\alpha\beta}\psi^I\partial_\alpha\partial_\beta\psi^I-\tau^2t^{-\delta-1}(2-t^{-\delta-1})(\psi^I)^2 \ dxdt\\
&=-\khw{P_1\psi^I,t^{-2}\psi^I}.}
Thus,
\fm{&\normw{r^{-1}\slashed\partial\psi^I}^2\lesssim \normw{t^{-1}\slashed\partial\psi^I}^2\\
&=-\khw{t^{-1}P_1\psi^I,t^{-1}\psi^I}-(4\tau^2)^{-1}\normw{t^{-1}P_2\psi^I}^2-\tau^{-1}\khw{P_2\psi^I,t^{-2}(1-t^{-\delta-1})\partial_t\psi^I}\\
&\quad+\int_{\Omega_{T_0,R''}}\tau^2 t^{-\delta-1}(2-t^{-\delta-1})(\psi^I)^2+t^{-\delta-1}(2-t^{-\delta-1})(\partial_t\psi^I)^2\ dxdt\\
&\lesssim \normw{t^{-1}P_1\psi^I}^2+\normw{t^{-1}\psi^I}^2+\tau^{-1}(\normw{P_2\psi^I}^2+\normw{t^{-2}\partial_t\psi^I}^2)\\
&\quad+\int_{\Omega_{T_0,R''}}\tau^2 t^{-\delta-1}(\psi^I)^2+t^{-\delta-1}(\partial_t\psi^I)^2\ dxdt\\
&\lesssim T_0^{-2}\normw{P_1\psi^I}^2+\normw{r^{-1}\psi^I}^2+\tau^{-1}(\normw{P_2\psi^I}^2+\normw{t^{-2}\partial_t\psi^I}^2)\\
&\quad+\tau^2\normw{t^{-\delta/2-3/2}\psi^I}^2+\normw{t^{-\delta/2-3/2}\partial_t\psi^I}^2.}
Apply Lemma \ref{lem3.3.1}. We conclude that for $\tau\gg1$,
\fm{&\normw{r^{-1}\slashed\partial\psi^I}^2\\&\lesssim T_0^{-2}\normw{P_1\psi^I}^2+\tau^{-1} \normw{P_2\psi^I}^2+\tau^2T_0^{-1}\normw{t^{-\delta/2-1}\psi^I}^2+(1+\tau^{-1}T_0^{\delta-1})T_0^{-1}\normw{t^{-\delta/2-1}\partial_t\psi^I}^2.}
We finish the proof by noticing that $T_0^{\delta-1}\leq 1$ and $\tau\gg 1$.
\end{proof}\rm
\bigskip

Let us now finish the proof of \eqref{prop3.7f1} and \eqref{prop3.7f2}. By \eqref{prop3.1a1} and  Lemma \ref{lem3.3.2}, we have \fm{\normw{X_\omega\psi^J}&\lesssim \normw{r^{-1}\slashed\partial\psi^J}\\
&\lesssim T_0^{-1}\normw{P_1\psi^J}+\tau^{-1/2} \normw{P_2\psi^J}+\tau T_0^{-1/2}\normw{t^{-\delta/2-1}\psi^J}+T_0^{-1/2}\normw{t^{-\delta/2-1}\partial_t\psi^J}.} By \eqref{prop3.1a1}, \eqref{prop3.1a2}, we have
\fm{&\normw{\frac{1}{2}V^{I,\alpha}_J\wt{\omega}_\alpha(\partial_t+\partial_r)\psi^J+\frac{1}{2}V^{I,\alpha}_J\widehat{\omega}_\alpha(-\partial_t+\partial_r)\psi^J}\\&\leq \frac{1}{2}C_1\normw{(1+t+r)^{-1}(\partial_t+\partial_r)\psi^J}+\frac{1}{2}C_1\normw{(1+t+r)^{-1-\gamma}(-\partial_t+\partial_r)\psi^J}\\
&\leq\frac{1}{2}C_1\normw{(1+t+r)^{-1}(-\frac{P_2\psi^J}{2\tau}+t^{-\delta-1}\partial_t\psi^J)}\\&\quad+\frac{1}{2}C_1\normw{(1+t+r)^{-1-\gamma}(-\frac{P_2\psi^J}{2\tau}-(2-t^{-\delta-1})\partial_t\psi^J)}\\
&\leq \frac{C_1}{4\tau}\cdot T_0^{-1}\normw{P_2\psi^J}+\frac{C_1}{2}\normw{t^{-\delta-2}\partial_t\psi^J}+C_1\normw{t^{-1-\gamma}\partial_t\psi^J}\\
&\leq C_1\tau^{-1}T_0^{-1}\normw{P_2\psi^J}+C_1T_0^{-1-\delta/2}\normw{t^{-1-\delta/2}\partial_t\psi^J}+C_1T_0^{\delta/2-\gamma}\normw{t^{-1-\delta/2}\partial_t\psi^J}.}
Since $\delta<\gamma$, we have $\delta/2-\gamma<-\gamma/2$. Then,
\fm{\normw{V^{I,\alpha}_J\partial_\alpha\psi^J}&\lesssim T_0^{-1}\normw{P_1\psi^J}+(\tau^{-1/2}+C_1\tau^{-1}T_0^{-1}) \normw{P_2\psi^J}+\tau T_0^{-1/2}\normw{t^{-\delta/2-1}\psi^J}\\&\quad+(C_1T_0^{-1-\delta/2}+C_1T_0^{-\gamma/2}+T_0^{-1/2})\normw{t^{-1-\delta/2}\partial_t\psi^J}.}
By choosing $T_0\gg_{C_1,R,R'',\gamma}1$, we have
\eq{\label{riestf2}\normw{V^{I,\alpha}_J\partial_\alpha\psi^J}&\leq T_0^{-1/2}\normw{P_1\psi^J}+2\tau^{-1/2} \normw{P_2\psi^J}+ T_0^{-\gamma/4}(\tau\normw{t^{-\delta/2-1}\psi^J}+\normw{t^{-1-\delta/2}\partial_t\psi^J}).}
By applying \eqref{riestf1} and \eqref{riestf2}, and by choosing $T_0\gg_{C_1,R,R'',\gamma}1$ and $\tau\gg 1$, we conclude  \eqref{prop3.7f1}. By the triangle inequality and Lemma \ref{lem3.3.1}, we have
\fm{&\normw{P_2(\psi^I)+R^I(\psi)}^2\\&\geq |\normw{P_2(\psi^I)-2\tau r^{-1}\psi^I}-\normw{R^I(\psi)+2\tau r^{-1}\psi^I}|^2\\
&=\normw{P_2(\psi^I)-2\tau r^{-1}\psi^I}^2+\normw{R^I(\psi)+2\tau r^{-1}\psi^I}^2\\&\quad-2\normw{P_2(\psi^I)-2\tau r^{-1}\psi^I}\normw{R^I(\psi)+2\tau r^{-1}\psi^I}\\
&\geq |\normw{P_2(\psi^I)}-2\tau \normw{r^{-1}\psi^I}|^2+\normw{R^I(\psi)+2\tau r^{-1}\psi^I}^2\\&\quad-2(\normw{P_2(\psi^I)}+2\tau \normw{r^{-1}\psi^I})\normw{R^I(\psi)+2\tau r^{-1}\psi^I}\\
&\geq \frac{1}{25}\normw{P_2\psi^I}^2+\normw{R^I(\psi)+2\tau r^{-1}\psi^I}^2-\frac{18}{5}\normw{P_2\psi^I}\normw{R^I(\psi)+2\tau r^{-1}\psi^I}\\
&\geq  \frac{1}{100}\normw{P_2\psi^I}^2-C\normw{R^I(\psi)+2\tau r^{-1}\psi^I}^2.}
By \eqref{prop3.7f1}, for $T_0\gg_{R,R'',C_1,\gamma}1$ and $\tau\gg1$, we have
\fm{&\normw{P_2(\psi^I)+R^I(\psi)}^2\\
&\geq  \frac{1}{100}\normw{P_2\psi^I}^2-CT_0^{-\gamma/2}\sum_J(\normw{P_1\psi^J}^2+\normw{P_2\psi^J}^2+\tau\normw{t^{-\delta/2-1}\psi^J}^2+\normw{t^{-\delta/2-1}\partial_t\psi^J}^2)\\&\quad-C\tau^{-1}\sum_J\normw{P_2\psi^J}^2\\
&\geq \frac{1}{100}\normw{P_2\psi^I}^2-\frac{1}{1000N}\sum_{J}(\normw{P_1\psi^J}^2+\normw{P_2\psi^J}^2+\tau\normw{t^{-\delta/2-1}\psi^J}^2+\normw{t^{-\delta/2-1}\partial_t\psi^J}^2).}
 We thus obtain \eqref{prop3.7f2}.

\subsection{Proof for \eqref{prop3.1con}}\label{sec3.4}
Let us finish the proof for Proposition \ref{prop3.1}. We combine \eqref{allest}, \eqref{s3imp}, Proposition \ref{prop3.1.1} and Proposition \ref{prop3.7}. For simplicity, if $V=(V^I)$ is an $\R^N$-valued function, then  we set $\normw{V}^2=\sum_I\normw{V^I}^2$. So far we have proved
\fm{&\sum_{I=1}^N\normw{(\wt{\Box}_{\tau,\delta,R'}\psi)^I}^2\\&\geq\frac{999}{1000}\normw{P_1\psi}^2+(\frac{9}{1000}-\tau^{-1})\normw{P_2\psi}^2+(\frac{1}{4}\tau-\frac{1}{1000})(\normw{t^{-1-\delta/2}\partial_t\psi}^2+\tau^2\normw{t^{-1-\delta/2}\psi}^2)\\&\quad+\sum_I [4\khw{P_2\psi^I,(t^{-1}-r^{-1}+r^{-1}t^{-\delta-1})\partial_t\psi^I}+\khw{P_2(\psi^I),2t^{-2}\psi^I}\\&\qquad\qquad+\khw{P_1(\psi^I),2\tau( (\delta+1)t^{-\delta-2}+2 t^{-1}(1-t^{-\delta-1}))\psi^I+R^I(\psi)}].}
To continue, we note that 
\eq{\label{temf1}-4\khw{P_2\psi^I,(t^{-1}-r^{-1}+r^{-1}t^{-\delta-1})\partial_t\psi^I}
&\leq 4\normw{P_2\psi^I}\normw{(tr)^{-1}(r-t+t^{-\delta})\partial_t\psi^I}\\&\leq CT_0^{\delta/2-1}\normw{P_2\psi^I}\normw{t^{-1-\delta/2}\partial_t\psi^I}\\&\leq CT_0^{-1/2}(\normw{P_2\psi^I}^2+\normw{t^{-1-\delta/2}\partial_t\psi^I}^2),}
\eq{\label{temf2}-\khw{P_2(\psi^I),2t^{-2}\psi^I}&\leq 2\normw{P_2\psi^I}\normw{t^{-2}\psi^I}\leq CT_0^{\delta/2-1}\normw{P_2\psi^I}\normw{t^{-1-\delta/2}\psi^I}\\
&\leq CT_0^{-1/2}(\normw{P_2\psi^I}^2+\normw{t^{-1-\delta/2}\psi^I}^2),}
and
\eq{\label{temf3}&-\khw{P_1(\psi^I),2\tau( (\delta+1)t^{-\delta-2}+2 t^{-1}(1-t^{-\delta-1}))\psi^I+R^I(\psi)}\\
&=-\khw{P_1(\psi^I),2\tau( (\delta-1)t^{-\delta-2}+2 (tr)^{-1}(r-t))\psi^I+R^I(\psi)+2\tau r^{-1}\psi^I}\\
&\leq C\normw{P_1\psi^I}(\tau\normw{t^{-2}\psi^I}+\normw{R^I(\psi)+2\tau r^{-1}\psi^I})\\
&\leq \frac{1}{2}\normw{P_1\psi^I}^2+C(\tau^2T_0^{\delta-2}\normw{t^{-1-\delta/2}\psi^I}^2+\normw{R^I(\psi)+2\tau r^{-1}\psi^I}^2)\\
&\leq \frac{1}{2}\normw{P_1\psi^I}^2+C\tau^2T_0^{-1}\normw{t^{-1-\delta/2}\psi^I}^2\\&\quad+CT_0^{-\gamma/2}(\normw{P_1\psi}^2+\normw{P_2\psi}^2+\tau^2\normw{t^{-\delta/2-1}\psi}^2+\normw{t^{-\delta/2-1}\partial_t\psi}^2)+C\tau^{-1}\normw{P_2\psi}^2.}
Take the sum of \eqref{temf1}, \eqref{temf2} and \eqref{temf3} over all $1\leq I\leq N$. For $T_0\gg1$ and $\tau\gg1$, we conclude that this sum has a upper bound \fm{\frac{1}{2}\norm{P_1\psi}^2+\frac{1}{1000}(\normw{P_1\psi}^2+\normw{P_2\psi}^2)+(\tau^2\normw{t^{-\delta/2-1}\psi}^2+\frac{1}{1000}\normw{t^{-\delta/2-1}\partial_t\psi}^2).}
As a result,\fm{&\sum_{I=1}^N\normw{(\wt{\Box}_{\tau,\delta,R'}\psi)^I}^2\\&\geq\frac{249}{500}\normw{P_1\psi}^2+(\frac{1}{125}-\tau^{-1})\normw{P_2\psi}^2+(\frac{1}{4}\tau-\frac{1}{500})(\normw{t^{-1-\delta/2}\partial_t\psi}^2+\tau^2\normw{t^{-1-\delta/2}\psi}^2).}
By choosing $\tau\gg1$, we conclude that 
\fm{\normw{P_1\psi}^2+\normw{P_2\psi}^2+\tau(\normw{t^{-1-\delta/2}\partial_t\psi}^2+\tau^2\normw{t^{-1-\delta/2}\psi}^2)\lesssim \sum_{I=1}^N\normw{(\wt{\Box}_{\tau,\delta,R'}\psi)^I}^2.}This estimate is in fact stronger than \eqref{prop3.1con}.

\section{Application of the Carleman estimates}\label{sec4}
In this section, we seek to prove the following proposition.
\prop{\label{prop4.1}Fix $R,R''\in\R$ such that $R''<R$.  Also fix $0<\gamma,\gamma'\leq 1$ and $T_s>|R''|+1$. Suppose that we have continuous functions $V_{J}^{I,\alpha}$'s and $W^{I}_J$'s for $I,J=1,\dots,N$ and $\alpha=0,1,2,3$ defined in $\Omega_{T_s,R''}$, and that they satisfy the following pointwise estimates:
\eq{\label{prop4.1a1}\sum_{I,J=1}^N\sum_{\alpha=0}^3|V^{I,\alpha}_J|\leq C_1 (1+t+r)^{-1}\qquad\forall(t,x)\in\Omega_{T_s,R''},}
and
\eq{\label{prop4.1a2}\sum_{I,J=1}^N|\sum_{\alpha=0}^3V^{I,\alpha}_J\widehat{\omega}_\alpha|+\sum_{I,J=1}^N|W^{I}_J|\leq  C_1(1+t+r)^{-1-\gamma}\qquad\forall(t,x)\in\Omega_{T_s,R''}.}
Suppose that $\phi=(\phi^I)\in C^2(\Omega_{T_s,R''};\R^N)$ is a solution to the PDE
\eq{\label{prop4.1pde} (\wt{\Box}\phi)^I&=\Box \phi^I+\sum_{J=1}^N\sum_{\alpha=0}^3V_J^{I,\alpha} \partial_\alpha \phi^J+\sum_{J=1}^NW_J^I\phi^J=0,\qquad I=1,2,\dots,N.}
Assume that $\phi=0$ for all $r-t\geq R$, and that it satisfies the pointwise estimates:
\eq{\label{prop4.1a3}|\phi|+|\partial\phi|\leq C_2(1+t+r)^{-3/2-\gamma'}\qquad \forall(t,x)\in\Omega_{T_s,R''},}
and
\eq{\label{prop4.1a4}|(\partial_t+\partial_r)\phi|\leq C_2(1+t+r)^{-5/2-\gamma'}\qquad \forall(t,x)\in\Omega_{T_s,R''}.}

Then, there exist  a large time $T_0\gg_{T_s,R,R'',C_1,C_2,\gamma,\gamma'}1$ and a sufficiently small $0<\delta\ll_{T_0,R,R''}1$, such that $\phi(t,x)= 0$ for all $(t,x)\in\Omega_{T_0,R''}\cap\{f_\delta\geq R''\}$.
}\rm\bigskip

The proof of Proposition \ref{prop4.1} relies on the Carleman estimates proved in Section \ref{s3}. However, $\phi$ is not assumed to be compactly supported, so Proposition \ref{prop3.1} cannot be applied to $\phi$ directly. This is why we need to define a new function $\psi$ by multiplying $\phi$ by some cutoff functions. We will explain the details later in this section. In Section \ref{sec4.1}, we  define this function $\psi$. Note that $\psi$ is a function of two parameters $\tau$ and $T$. Here $\tau$ comes from the Carleman estimates, and $T$ is the time where we truncate the solution. In Section \ref{sec4.2}, we take $T\to\infty$ and then take $\tau\to\infty$ in the Carleman estimates. It turns out that the right hand side of \eqref{prop3.1con} has a limit $0$, which finishes the proof of Proposition \ref{prop4.1}.

\subsection{Setup}\label{sec4.1}
By Proposition \ref{prop3.1}, there exists a sufficiently large time $T_0\gg_{T_s,R,R'',C_1,\gamma}1$ such that we have the Carleman estimates \eqref{prop3.1con} for all $\tau\gg1$, $\delta\in(0,\gamma)$ and functions $\psi\in C_c^\infty(\Omega_{T_0,R''};\R^N)$. The choice of $T_0$ is independent of the choice of $\delta$, so we can choose $\delta$ after we obtain this time $T_0$. Because of the limit \fm{\lim_{\delta\to 0^+}\delta^{-1}T_0^{-\delta}=\infty,}we are able to choose $0<\delta\ll_{T_0,R,R''}1$ such that $\delta^{-1}T_0^{-\delta}>2(R-R'')$. We then define $f=f_\delta(t,x)$ by \eqref{fdefn} using this $\delta$. 

Fix any $R'\in\R$ such that $R''<R'<R$. Let $\chi_1,\chi_2\in C^\infty(\R)$ be two $\R$-valued cutoff functions satisfying the following properties $0\leq \chi_1(s),\chi_2(s)\leq 1$ for all $s\in\R$; $\chi_1(s)=0$ for all $s\leq (R''+R')/2$ and $\chi_1(s)=1$ for all $s\geq R'$; $\chi_2(s)=1$ for all $s\leq 1$ and $\chi_2(s)=0$ for all $s\geq 2$. Also fix a large time $T>T_0$. For all $\tau\gg 1$, we set
\eq{\label{s4psidefn} \psi:=\chi_1(f)\chi_2(t/T)e^{\tau(f-R')}\phi.}
From the choices of the cutoff functions, we have
$\psi=e^{\tau(f-R')}\phi$ whenever \fm{T_0\leq t\leq T,\qquad r-t-\delta^{-1}t^{-\delta}\geq R'.}
We also have $\psi=0$ whenever 
\fm{t\geq 2T,\qquad r-t-\delta^{-1}t^{-\delta}\leq (R''+R')/2.}
To visualize these cutoff regions, one could temporarily ignore the term $\delta^{-1}t^{-\delta}$ which is supposed to be negligible for sufficiently large time.

It is clear that $\psi\in C_c^2(\Omega_{T_0,R''};\R^N)$. In fact, we already know that $\phi=0$ for $r-t\geq R$, and that  $\psi=0$ whenever $t>2T$ or whenever $f(t,x)<(R''+R')/2$ because of the choice of $\chi_1,\chi_2$. Moreover, whenever $f(t,x)\geq R''$ and $T_0\leq t<\frac{101}{100}T_0$, we have \fm{r-t-(R-R'')>r-t-\delta^{-1}T_0^{-\delta}\cdot (\frac{100}{101})^{\delta}>r-t-\delta^{-1}t^{-\delta}>R''\Longrightarrow r\geq t+R.}
So the assumptions on $\delta$ above guarantee that $\phi\equiv 0$ in a small neighborhood of the surface $\{t=T_0,\ f(t,x)\geq R''\}$.  In summary, we have $\psi\in C_c^2(\Omega_{T_0,R''};\R^N)$.

\subsection{Proof of Proposition \ref{prop4.1}}\label{sec4.2} We can now apply the Carleman estimates. By \eqref{prop3.1con}, we conclude that for all $\tau\gg 1$, we have
\eq{\label{s4carest}\tau^3\int_{\Omega_{T_0,R''}}|\psi|^2t^{-\delta}\ dxdt\lesssim \int_{\Omega_{T_0,R''}}|\wt{\Box}_{\tau,\delta,R'}\psi|^2t^2\ dxdt.}
Note that the constant in this estimate is independent of $\tau$, $T_0$, $\delta$ and $\psi$. Here $\wt{\Box}_{\tau,\delta,R'}$ is defined by \eqref{boxfdefn}, so we have
\fm{\wt{\Box}_{\tau,\delta,R'}\psi&=e^{\tau(f-R')}\wt{\Box}(\chi_1(f)\chi_2(t/T)\phi).}
That is, for each $I=1,2,\dots,N$, we have
\eq{&(\wt{\Box}_{\tau,\delta,R'}\psi)^I\\
&=e^{\tau(f-R')}[\Box(\chi_1(f)\chi_2(t/T)\phi^I)+\sum_{J,\alpha}V_J^{I,\alpha}\partial_\alpha(\chi_1(f)\chi_2(t/T)\phi^J)+\sum_{J}W^I_J\chi_1(f)\chi_2(t/T)\phi^J]\\
&=e^{\tau(f-R')}\chi_1(f)\chi_2(t/T)(\wt{\Box}\phi)^I+\mcl{R}_1+\mcl{R}_2=\mcl{R}_1+\mcl{R}_2.}
The last identity follows because $\phi$ is a solution to $\wt{\Box}\phi=0$. Here we define 
\fm{\mcl{R}_1&:=e^{\tau(f-R')}[\chi_2(t/T)\phi^I\cdot\Box(\chi_1(f))+2\partial^\beta(\chi_1(f))\partial_\beta(\chi_2(t/T)\phi^I)+V_J^{I,\alpha}\partial_\alpha(\chi_1(f))\chi_2(t/T)\phi^J]}
and
\fm{\mcl{R}_2&:=e^{\tau(f-R')}\chi_1(f)[\Box(\chi_2(t/T))\phi^I+2\partial^\beta(\chi_2(t/T))(\partial_\beta\phi^I)+V_J^{I,\alpha}\partial_\alpha(\chi_2(t/T))\phi^J].}
Let us briefly explain why the error term here is written as the sum of $\mcl{R}_1$ and $\mcl{R}_2$. In fact, every term involving a derivative of $\chi_1(f)$ is included in $\mcl{R}_1$. Because of the definition of $\chi_1$, we have $\mcl{R}_1\neq0$ only if $f-R'<0$ and thus $\lim_{\tau\to\infty}e^{\tau(f-R')}=0$. All the other terms are put in $\mcl{R}_2$, and we notice that all such terms involve a derivative of $\chi_2(t/T)$. Because of the definition of $\chi_2$, we have $\mcl{R}_2\neq0$ only if $t\sim T$.

We have
\fm{\mcl{R}_1&=e^{\tau(f-R')}\chi_2(t/T)\phi^I(\chi_1'(f)\Box f+\chi_1''(f)\partial_\alpha f\partial^\alpha f)\\&\quad+2e^{\tau(f-R')}\chi_1'(f)\chi_2(t/T)\partial^\beta f\partial_\beta \phi^I-2e^{\tau(f-R')}\chi_1'(f)\partial_t f \chi_2'(t/T) T^{-1}\phi^I\\
&\quad+e^{\tau(f-R')}V_J^{I,\alpha}\partial_\alpha f\chi_1'(f)\chi_2(t/T)\phi^J\\
&=e^{\tau(f-R')}\chi_2(t/T)\phi^I(\chi_1'(f)(2r^{-1}+(\delta+1) t^{-\delta-2})+\chi_1''(f)t^{-\delta-1}(2-t^{-\delta-1}))\\&\quad+2e^{\tau(f-R')}\chi_1'(f)[\chi_2(t/T)((\partial_t+\partial_r) \phi^I-t^{-\delta-1}\partial_t\phi^I)+2T^{-1}(1-t^{-\delta-1}) \chi_2'(t/T) \phi^I]\\
&\quad+e^{\tau(f-R')}(V_J^{I,\alpha}\widehat{\omega}_\alpha+V_J^{I,0}t^{-\delta-1})\chi_1'(f)\chi_2(t/T)\phi^J}
and
\fm{\mcl{R}_2&=e^{\tau(f-R')}\chi_1(f)[-T^{-2}\chi_{2}''(t/T)\phi^I-2 T^{-1}\chi_2'(t/T) \partial_t\phi^I +V_J^{I,0}T^{-1}\chi_2'(t/T)\phi^J].}
Recall the estimates \eqref{prop4.1a1} and \eqref{prop4.1a2} for $V_{*}^{*,*}$ and $W^*_*$, and the estimates \eqref{prop4.1a3} and \eqref{prop4.1a4} for~$\phi$. Also recall that $t\sim r$ in $\Omega_{T_0,R''}$ and that $t\sim T$ in the support of $\chi_2'(t/T)$. As a result, we have
\fm{|\mcl{R}_1|\lesssim_{C_1,C_2} e^{\tau(f-R')}(|\chi_2(t/T)|+|\chi_2'(t/T)|)\cdot(|\chi_1'(f)|+|\chi_1''(f)|)\cdot t^{-5/2-\gamma'}1_{r-t\leq R},}
\fm{|\mcl{R}_2|&\lesssim_{C_1,C_2} e^{\tau(f-R')}|\chi_1(f)|(|\chi_{2}'(t/T)|+|\chi_{2}''(t/T)|)\cdot t^{-5/2-\gamma'}1_{r-t\leq R}.}
Thus,
\eq{&\normw{\mcl{R}_1}^2+\normw{\mcl{R}_2}^2\\&\lesssim \int_{\Omega_{T_0,R''}}e^{2\tau(f-R')}(|\chi_2(t/T)|+|\chi_2'(t/T)|)^2\cdot(|\chi_1'(f)|+|\chi_1''(f)|)^2t^{-3-2\gamma'}1_{r-t\leq R}\\&\qquad\qquad+e^{2\tau(f-R')}|\chi_1(f)|^2(|\chi_{2}'(t/T)|+|\chi_{2}''(t/T)|)^2t^{-3-2\gamma'}1_{r-t\leq R}\  dxdt\\
&\lesssim \int_{T_0}^{\infty}\int_{R''\leq r-t\leq R}e^{2\tau(f-R')}t^{-3-2\gamma'}1_{ f(t,x)<R'}\ dxdt+\int_{T}^{2T}\int_{R''\leq r-t\leq R}e^{2\tau(f-R')}t^{-3-2\gamma'} dxdt.}
Note that 
\fm{&\int_{T_0}^{\infty}\int_{R''\leq r-t\leq R}e^{2\tau(f-R')}t^{-3-2\gamma'}1_{ f(t,x)<R'}\ dxdt\\&\lesssim\int_{T_0}^{\infty} \int_{R''\leq r-t\leq R}t^{-3-2\gamma'}\ dxdt\lesssim_{R,R''}\int_{T_0}^{\infty} t^2\cdot t^{-3-2\gamma'}\ dt\lesssim_{R,R'',T_0,\gamma'}1, }
\fm{&\int_{T}^{2T}\int_{R''\leq r-t\leq R}e^{2\tau(f-R')}t^{-3-2\gamma'}\ dxdt\\&\lesssim e^{\tau(|R|+|R''|+|R'|)}\int_{T}^{2T}\int_{R''\leq r-t\leq R}t^{-3-2\gamma'}\ dxdt\lesssim_{R,R''} e^{\tau(|R|+|R''|+|R'|)}\int_{T}^{2T}t^{2}\cdot t^{-3-2\gamma'}\ dt\\
&\lesssim_{\gamma'} e^{\tau(|R|+|R''|+|R'|)} T^{-2\gamma'}.}
By the Lebesgue dominated convergence theorem, we have
\fm{&\lim_{\tau\to\infty}\int_{T_0}^{\infty}\int_{R''\leq r-t\leq R}e^{2\tau(f-R')}t^{-3-2\gamma'}1_{ f(t,x)<R'}\ dxdt\\&=\int_{T_0}^{\infty}\int_{R''\leq r-t\leq R}\lim_{\tau\to\infty}e^{2\tau(f-R')}t^{-3-2\gamma'}1_{ f(t,x)<R'}\ dxdt=0,}
\fm{\lim_{T\to\infty}\int_{T}^{2T}\int_{R''\leq r-t\leq R}e^{2\tau(f-R')}t^{-3-2\gamma'}\ dxdt=0.} 
In conclusion, we have
\fm{\lim_{\tau\to\infty}\limsup_{T\to\infty}(\normw{\mcl{R}_1}^2+\normw{\mcl{R}_2}^2)=0}
and thus
\eq{\lim_{\tau\to\infty}\limsup_{T\to\infty}\normw{\wt{\Box}_{\tau,\delta,R'}\psi}^2=0.}
By \eqref{s4carest} we deduce that
\eq{\label{limcar}\lim_{\tau\to\infty}\limsup_{T\to\infty}\tau^3\int_{\Omega_{T_0,R''}}|\chi_1(f)\chi_2(t/T)e^{\tau(f-R')}\phi|^2t^{-\delta}\ dxdt=0.}
By the Fatou's lemma, we have
\fm{&\int_{\Omega_{T_0,R''}}|\chi_1(f)e^{\tau(f-R')}\phi|^2t^{-\delta}\ dxdt\\
&=\int_{\Omega_{T_0,R''}}\lim_{T\to\infty}|\chi_1(f)\chi_2(t/T)e^{\tau(f-R')}\phi|^2t^{-\delta}\ dxdt\\
&\leq \liminf_{T\to\infty}\int_{\Omega_{T_0,R''}}|\chi_1(f)\chi_2(t/T)e^{\tau(f-R')}\phi|^2t^{-\delta}\ dxdt\\
&\leq \limsup_{T\to\infty}\int_{\Omega_{T_0,R''}}|\chi_1(f)\chi_2(t/T)e^{\tau(f-R')}\phi|^2t^{-\delta}\ dxdt.}
It then follows from \eqref{limcar} that
\fm{\lim_{\tau\to\infty}\tau^3\int_{\Omega_{T_0,R''}}|\chi_1(f)e^{\tau(f-R')}\phi|^2t^{-\delta}\ dxdt=0.}
By the Fatou's lemma again, we conclude that
\eq{\label{s4estlim} \int_{\Omega_{T_0,R''}}\liminf_{\tau\to\infty}\tau^3|\chi_1(f)e^{\tau(f-R')}\phi|^2t^{-\delta}\ dxdt=0.}
However, this limit forces $\phi=0$ for all $(t,x)\in\Omega_{T_0,R''}$ such that $f(t,x)\geq R'$. Otherwise, we have $|\phi|>0$ in a nonempty open set $U$ in $\R^{1+3}$ by continuity. In this case, the integrand on the left side of \eqref{s4estlim} is infinite on a set of positive measure. This contradicts with the limit \eqref{s4estlim}.

Finally, we notice that $R''<R'<R$ can be chosen arbitrarily. Thus $\phi\equiv 0$ in the domain \fm{\bigcup_{R'\in(R'',R)}\{(t,x)\in\Omega_{T_0,R''}:\ f(t,x)\geq R'\}=\{(t,x)\in\Omega_{T_0,R''}:\ f(t,x)>R''\}.}
By continuity of $\phi$, we also have $\phi(t,x)= 0$ for $(t,x)\in\Omega_{T_0,R''}$ such that $f(t,x)=R''$. This finishes the proof of Proposition \ref{prop4.1}. 

\section{Proof of the main theorems}\label{sec5}

In this section we finish the proof of Theorem \ref{mthm0} and Theorem \ref{mthm}.  In Section \ref{sec5.2} and Section \ref{sec5.3}, we prove Theorem \ref{mthm} and Theorem \ref{mthm0}, respectively. At the end of Section \ref{sec5.3}, we also prove the results stated in Remark \ref{rmkk1.4}.

\subsection{Vanishing in $\D$}\label{sec5.2} In this section we finish the proof of Theorem \ref{mthm}. Fix two constants $R_1,R_2\in\R$ such that $R_1>0$ and $|R_2|<R_1$. Recall from Theorem \ref{mthm} that we  define the open set
\eq{\label{ddefn} \D=\D_{R_1,R_2}:=\{(t,x)\in\R^{1+3}:\ t> 0,\ (r-t-R_2)(r+t+R_2)>R_1^2-R_2^2\}.}
Similarly, for each $R''\in(R_2,R_1)$, we define $\D_{R_1,R''}$  by \eqref{ddefn} with $R_2$ replaced by $R''$. Note that $\bigcup_{R''\in (R_2,R_1)}\D_{R_1,R''}=\D_{R_1,R_2}$. In fact, $(t,x)\in\D_{R_1,R_2}$ if and only if $t>0$ and
\fm{R_2<\frac{r^2-R_1^2-t^2}{2t}.}
For some sufficiently small $\eps>0$, we have
\fm{R_2+\eps<\frac{r^2-R_1^2-t^2}{2t},}which implies that $(t,x)\in\D_{R_1,R_2+\eps}$.

Let $\phi$ be a solution to $\wt{\Box}\phi=0$ with the properties stated in the Theorem \ref{mthm}.  Our goal  is to show $\phi\equiv 0$ in $\D$. Because of the continuity of $\phi$, it suffices to prove  the following proposition.

\prop{\label{prop5.6} For each $R''\in(R_2,R_1)$, we have $\phi\equiv 0$ in $\D_{R_1,R''}$. }\rm\bigskip

In the rest of this section, we  will prove this proposition. Now we fix  $R''\in(R_2,R_1)$ and   set  \fm{T_s:=\frac{R_1^2-(R'')^2}{2(R''-R_2)}>0.}
We claim that $\Omega_{T_s,R''}\subset\D$. In fact, for each $(t,x)\in\Omega_{T_s,R''}$, we have
\fm{(r-t-R_2)(r+t+R_2)&> (R''-R_2)(R''+2T_s+R_2)=(R'')^2-R_2^2+2T_s(R''-R_2)\\
&=R_1^2-R_2^2.}
That is, we have $(t,x)\in\D$. Given such a pair of $(T_s,R'')$, we apply Proposition \ref{prop4.1}. Our conclusion is that  $\phi\equiv 0$ in $\Omega_{T_0,R''}\cap\{f_\delta\geq R''\}$ for some $T_0\gg_{T_s,R_1,R'',C_1,C_2,\gamma,\gamma'}1$ and $0<\delta\ll_{T_0,R_1,R''}1$. In summary, so far  we have proved that $\phi\equiv 0$ in 
\eq{\label{ff510}K_{\delta,T_0,R'',R_1}:=(\Omega_{T_0,R''}\cap\{f_\delta\geq R''\})\cup\{(t,x)\in\R^{1+3}:\ t>0,\ r-t\geq R_1\}.}

To continue, for fixed constants $\nu>0$ and $\kappa\in\R$, we define a family of surfaces 
\eq{\label{f511}S_{\nu,\kappa,c}:=\{(t,x)\in\R^{1+3}:\ t>0,\ (r+\nu)^2-(t+\kappa)^2=c\},\qquad c>0.}
According to Example \ref{exmhyp}, these surfaces are strongly pseudoconvex in $\bigcup_{c>0}S_{\nu,\kappa,c}$. We seek to foliate the  region $\D_{R_1,R''}$ with this family of surfaces $\{S_{\nu,\kappa,c}\}$ where $(\nu,\kappa,c)$ satsifies some constraints.

We first specify the constraints on $(\nu,\kappa,c)$. 

\lem{\label{lem5.7}We have \eq{\label{lem5.7f1}\D_{R_1,R''}=\bigcup\{S_{\nu,\kappa,c}:\ \nu>0,\ R''+\nu<\kappa<R_1+\nu,\ c>(R_1+\nu)^2-\kappa^2\}.}
Here we note that $(R_1+\nu)^2-\kappa^2>0$.}
\begin{proof} We first explain why $\nu>0$ and $R''+\nu<\kappa<R_1+\nu$ implies $(R_1+\nu)^2-\kappa^2>0$. In fact,  we have $|R''|<R_1$ since $R_1>|R_2|$ and $R_1>R''>R_2$. Thus, if $R''+\nu<\kappa<R_1+\nu$, then either $0\leq \kappa<R_1+\nu$ or $0>\kappa>R''+\nu>-R_1+\nu$. In summary, we must have $|\kappa|<R_1+\nu$.

Fix $\nu>0$, $R''+\nu<\kappa<R_1+\nu$ and $c>(R_1+\nu)^2-\kappa^2$. We claim that $S_{\nu,\kappa,\delta}\subset \D_{R_1,R''}$. To prove this claim, we fix $(t,x)\in S_{\nu,\kappa,c}$. Since $t>0$ and $|R''|<R_1$, we have
\fm{(r+t+R'')(r-t-R'')>(R_1+R'')(R_1-R'')=R_1^2-(R'')^2} whenever $r-t> R_1$. If $r-t\leq R_1$, we have \fm{& (r+t+R'')(r-t-R'')-[(r+\nu)^2-(t+\kappa)^2]\\&=-2r\nu+2t(\kappa-R'')-\nu^2-(R'')^2+\kappa^2\\
&\geq-2(r-t)\nu-\nu^2-(R'')^2+\kappa^2\geq-2R_1\nu+\nu^2+(R'')^2-\kappa^2.}
As a result, for each $(t,x)\in S_{\nu,\kappa,c}\cap\{r-t\leq R_1\}$, we have
\fm{(r+t+R'')(r-t-R'')\geq c-2R_1\nu-\nu^2-(R'')^2+\kappa^2.}
If $c>(R_1+\nu)^2-\kappa^2$, the right side of this inequality is larger than $R_1^2-(R'')^2$. We thus conclude that $S_{\nu,\kappa,c}\subset \D_{R_1,R''}$.

Conversely, we fix $(t,x)\in\D_{R_1,R''}$. Since $(r+t+R'')(r-t-R'')>R_1^2-R_2^2$, we can choose $\eps>0$ so that $(r+t+R'')(r-t-R'')>R_1^2-R_2^2+\eps$. We emphasize that $\eps$ is chosen before $(\nu,\kappa,c)$ is chosen. According to the computations above, we have
\fm{&(r+\nu)^2-(t+\kappa)^2\\&=(r+t+R'')(r-t-R'')+(R'')^2+2(r-t)\nu+2t(\nu+R''-\kappa)+\nu^2-\kappa^2\\
&>R_1^2+\eps+2(r-t)\nu+2t(\nu+R''-\kappa)+\nu^2-\kappa^2\\
&=(R_1+\nu)^2-\kappa^2+\eps+2(r-t-R_1)\nu+2t(\nu+R''-\kappa).}
We can choose $\nu>0$ and $R''+\nu<\kappa<R_1+\nu$ (both $\nu$ and $\kappa$ depend on $(t,x),\ \eps$ and $R''$) such that $\eps+2(r-t-R_1)\nu+2t(\nu+R''-\kappa)>0$. As a result, we have $(t,x)\in S_{\nu,\kappa,c}$ with $c>(R_1+\nu)^2-\kappa^2>0$.

\end{proof}\rm\bigskip

We hope to apply Corollary \ref{cors5carest} on these $S_{\nu,\kappa,c}$. However, Corollary \ref{cors5carest} is a local result, so some type of compactness would be necessary.

\lem{\label{lem5.8} Let $\nu,c>0$ and $\kappa\in\R$ be constants such that $R''+\nu<\kappa<R_1+\nu$. Then, for any fixed two constants $c_1,c_2$ such that $(R_1+\nu)^2-\kappa^2<c_1\leq c_2$, there exist two constants $0<T_1<T_2$, such that  \fm{\bigcup_{c_1\leq c\leq c_2}S_{\nu,\kappa,c}\setminus K_{\delta,T_0,R'',R_1}\subset[T_1,T_2]\times \R^3.} As a result, the closure of $\bigcup_{c_1\leq c\leq c_2}S_{\nu,\kappa,c}\setminus K_{\delta,T_0,R'',R_1}$ is a compact subset of $[T_1,T_2]\times \R^3$.}
\begin{proof} Let $(t,x)$ be any point in $\bigcup_{c_1\leq c\leq c_2}S_{\nu,\kappa,c}\setminus K_{\delta,T_0,R'',R_1}$. We first show that $t\geq T_1>0$ for some fixed time $T_1>0$. If this is false,  we can find a sequence $\{(t_n,x_n)\}$ such that $(t_n,x_n)\in \bigcup_{c_1\leq c\leq c_2}S_{\nu,\kappa,c}\setminus K_{\delta,T_0,R'',R_1}$ for all $n$ and $\lim_{n\to\infty}t_n=0$. It follows from \eqref{f511} that
\fm{c_1\leq (|x_n|+\nu)^2-(t_n+\kappa)^2\leq c_2,\qquad \forall n=1,2,\dots.}
It is clear that $\{x_n\}$ is a bounded sequence in $\R^3$, so it has a convergent subsequence. Without loss of generality, we assume that $\{x_n\}$ converges to $y\in\R^3$. As a result,
\fm{c_1\leq (|y|+\nu)^2-\kappa^2\leq c_2.}
Since $\nu>0$ and  $(R_1+\nu)^2-\kappa^2<c_1$, we have $|y|>R_1$. It follows that  $\lim_{n\to\infty}(|x_n|-t_n)=|y|>R_1$. But this implies that $|x_n|-t_n>R_1$ for all sufficiently large $n$. That is, $(t_n,x_n)\in K_{\delta,T_0,R'',R_1}$ for all sufficiently large $n$. A contradiction. As a result, we have $\bigcup_{c_1\leq c\leq c_2}S_{\nu,\kappa,c}\setminus K_{\delta,T_0,R'',R_1}\subset [T_1,\infty)\times\R^3$ for some $T_1>0$.

Next we show that $t\leq T_2$ for some fixed time $T_2\in(0,\infty)$. If this is false,  we can find a sequence $\{(t_n,x_n)\}$ such that $(t_n,x_n)\in \bigcup_{c_1\leq c\leq c_2}S_{\nu,\kappa,c}\setminus K_{\delta,T_0,R'',R_1}$ for all $n$ and $\lim_{n\to\infty}t_n=\infty$. Without loss of generality, we assume that $t_n>T_0$ for all $n$. Again, we have
\fm{c_1\leq (|x_n|+t_n+\kappa+\nu)(|x_n|-t_n+\nu-\kappa)=(|x_n|+\nu)^2-(t+\kappa)^2\leq c_2,\qquad \forall n=1,2,\dots.}
Note that this equation implies that $\lim_{n\to\infty}|x_n|=\infty$.
In addition,  we should have $f_\delta(t_n,x_n)<R''$ for all $n$. That is,
\fm{|x_n|-t_n-\delta^{-1}t_n^{-\delta}<R'',\qquad \forall n=1,2,\dots.}
As a result, 
\fm{\frac{c_2}{|x_n|+t_n+\kappa-\nu}+\kappa-\nu-\delta^{-1}t_n^{-\delta}<R''.}
By sending $n\to\infty$, we conclude that $\kappa-\nu\leq R''$. This contradicts with our assumptions. As a result, we have $\bigcup_{c_1\leq c\leq c_2}S_{\nu,\kappa,c}\setminus K_{\delta,T_0,R'',R_1}\subset [T_1,T_2]\times\R^3$ for some $T_2>0$.
\end{proof}\rm\bigskip

In the next lemma, we show that $S_{\nu,\kappa,c}\subset K_{\delta,T_0,R'',R_1}$ for sufficiently large $c>0$.

\lem{\label{lem5.9}Let $\nu>0$ and $\kappa\in\R$ be constants such that $R''+\nu<\kappa<R_1+\nu$. Then, for sufficiently large $c\gg_{\delta,T_0,R_1,\nu,\kappa}1$ (in particular, $c>(R_1+\nu)^2-\kappa^2>0$), we have $S_{\nu,\kappa,c}\subset K_{\delta,T_0,R'',R_1}$.} 
\begin{proof}First, on $S_{\nu,\kappa,c}$ we have
\fm{(r-t+\nu-\kappa)(r+t+\nu+\kappa)=c.}
Here $c>0$. If  $r-t+\nu-\kappa<0$ and $r+t+\nu+\kappa<0$, then  $0<2r+2\nu<0$ which is impossible. As a result, we have $r-t+\nu-\kappa>0$ and $r+t+\nu+\kappa>0$. 

Now fix $(t,x)\in S_{\nu,\kappa,c}$. If $r-t\geq R_1$, we conclude that $(t,x)\in K_{\delta,T_0,R'',R_1}$ by \eqref{ff510}. So let us assume that $r-t<R_1$. That is,
\fm{r-t=\frac{c}{r+t+\nu+\kappa}+\kappa-\nu<R_1.}
Since $r+t+\nu+\kappa>0$ and $R_1+\nu>\kappa$, it follows that
\fm{\frac{c}{R_1+\nu-\kappa}<r+t+\nu+\kappa<2t+\nu+\kappa+R_1.}
Moreover, we have
\fm{f_\delta(t,x)&=r-t-\delta^{-1}t^{-\delta}=\frac{c}{r+t+\nu+\kappa}+\kappa-\nu-\delta^{-1}t^{-\delta}>\kappa-\nu-\delta^{-1}t^{-\delta}.}
By choosing $c\gg_{\delta,T_0,R_1,\nu,\kappa}1$, we have $t>T_0$ and $0<\delta^{-1}t^{-\delta}<\kappa-\nu-R''$. This constant $c$ does not depend on $(t,x)$. Thus, $(t,x)\in \Omega_{T_0,R''}\cap\{f_\delta\geq R''\}$. This finishes the proof.

\end{proof}\rm\bigskip

Let us prove Proposition \ref{prop5.6}. Fix $\nu>0$ and  $R''+\nu<\kappa<R_1+\nu$. Set
\eq{I:=\{c_0>(R_1+\nu)^2-\kappa^2:\ \phi=0\quad \text{on }\bigcup_{c\geq c_0}S_{\nu,\kappa,c}\}.}
By Lemma \ref{lem5.9}, we have $I\neq \varnothing$, so we can set $c_0:=\inf I\geq (R_1+\nu)^2-\kappa^2$. We claim that  $c_0=(R_1+\nu)^2-\kappa^2$.

Let us instead assume that $c_0>(R_1+\nu)^2-\kappa^2$. For any $(R_1+\nu)^2-\kappa^2<c'<c_0$, we have $c'\notin I$, so there exists $c'\leq c< c_0$ such that $\{\phi\neq 0\}\cap S_{\nu,\kappa,c}\neq\varnothing$. We thus obtain a sequence of points $\{(t_n,x_n)\}$ and a sequence of real numbers $\{c_n\}$, such that $(t_n,x_n)\in S_{\nu,\kappa,c_n}$, $\phi(t_n,x_n)\neq 0$, $(R_1+\nu)^2-\kappa^2<c_n< c_0$ and $\lim_{n\to\infty}c_n=c_0$.  Here we can choose $c_2'\geq c_1'>(R_1+\nu)^2-\kappa^2$ such that $c_1'\leq c_n\leq c_2'$ for each $n$. And since $\phi\equiv 0$ on $K_{\delta,T_0,R'',R_1}$, for each $n$ we have $(t_n,x_n)\in \bigcup_{c_1'\leq c\leq c_2'}S_{\nu,\kappa,c}\setminus K_{\delta,T_0,R'',R_1}$. By Lemma \ref{lem5.8}, the closure of $\bigcup_{c_1'\leq c\leq c_2'}S_{\nu,\kappa,c}\setminus K_{\delta,T_0,R'',R_1}$ is a compact set contained in $[T_1,T_2]\times \R^3$ for some $0<T_1<T_2<\infty$. Thus, a  subsequence of $\{(t_n,x_n)\}$ converges. Without loss of generality, we assume that the sequence $\{(t_n,x_n)\}$ itself converges to $(t_\infty,x_\infty)$.

Now, it is clear that $(t_\infty,x_\infty)\in([T_1,T_2]\times\R^3)\cap S_{\nu,\kappa,c_0}$, so there exists an open neighborhood of $(t_\infty,x_\infty)$ contained in \fm{\mcl{O}=\{(t,x)\in\R^{1+3}:\ t>0,\ (r+\nu)^2-(t+\kappa)^2>0\}.} Since $\phi=0$ on $\bigcup_{c\geq c_0}S_{\nu,\kappa,c}$, we can apply Corollary \ref{cors5carest} at $(t_\infty,x_\infty)$. As a result, there exists an open neighborhood $\mcl{N}$ of $(t_\infty,x_\infty)$, such that $\phi\equiv 0$ in $\mcl{N}$. However, since the sequence $\{(t_n,x_n)\}$ converges to $(t_\infty,x_\infty)$, so $(t_n,x_n)\in \mcl{N}$ for all sufficiently large $n$. We obtain a contradiction as $\phi(t_n,x_n)\neq 0$.

So far, we have proved that $\inf I=(R_1+\nu)^2+\kappa^2$, so $\phi=0$ on $\bigcup_{c>(R_1+\nu)^2+\kappa^2}S_{\nu,\kappa,c}$. Since $\nu>0$ and $R''+\nu<\kappa<R_1+\nu$ can be chosen arbitrarily, we conclude that $\phi\equiv 0$ in $\D_{R_1,R''}$ by applying Lemma \ref{lem5.7}. This ends the proof of Proposition \ref{prop5.6}.

\subsection{Proof of Theorem \ref{mthm0} and Remark \ref{rmkk1.4}}\label{sec5.3}

In this section we explain why Theorem \ref{mthm} implies Theorem \ref{mthm0}. We also present a brief proof of the result stated in Remark \ref{rmkk1.4}.

\subsubsection{Proof of part (i) in Theorem \ref{mthm0}}
Let $u$ and $\wt{u}$ be two smooth global solutions as in Theorem \ref{mthm0}. Let $R_1$ and $R_2$ be the corresponding constants as in Theorem \ref{mthm0}, and we  assume that $|R_2|<R_1$. Set $\phi=u-\wt{u}$. It follows that for each $I=1,\dots,N$, 
\fm{\Box\phi^I&=\Box u^I-\Box\wt{u}^I=Q^I(u,\partial u)-Q^I(\wt{u},\partial\wt{u}).}
To continue, we recall a useful lemma. We remark that it can be viewed as a variant of  Theorem 1.1.9 in H\"{o}rmander \cite{hormI}.
\lem{\label{lem5.10horm}Fix two integers $K,M>0$. Suppose that $f=f(X)$ is a $C^M$ function near the origin in $\R^K$ and that $f(0)=0$. Then, for each $X,Y\in\R^K$ near the origin, we have \eq{f(X)-f(Y)=\sum_{i=1}^K (X_i-Y_i)f_{i}(X,Y).}
Here the $f_i$'s are some $C^{M-1}$ functions defined by
\eq{f_i(X,Y):=\int_0^1(\partial_{i}f)(\rho X+(1-\rho)Y)\ d\rho.}
 }\rm\bigskip

We recall that $Q^I=Q^I(u,v)$ is a function of \eq{(u,v)=((u^J)_{J=1,\dots, N},(v^{J}_\alpha)_{J=1,\dots,N;\ \alpha=0,1,2,3})\in\R^N\times\R^{4N}}
where the variable $u^J$ corresponds to the function $u^J(t,x)$ and the variable $v_\alpha^J$ corresponds to the derivative $(\partial_\alpha u^J)(t,x)$ in \eqref{swe}. By applying this lemma, we have
\eq{Q^I(\wt{u},\partial \wt{u})-Q^I(u,\partial u)&=\sum_{J=1}^NW_J^I\phi^J+\sum_{\alpha=0}^3\sum_{J=1}^N V^{I,\alpha}_J\partial_\alpha \phi^J.}
Here 
\eq{\label{sec5.3.1f3}W_J^I(t,x)&:=\int_0^1 (\partial_{u^J}Q^I)(\rho \wt{u}+(1-\rho) u,\rho\partial \wt{u}+(1-\rho)\partial u)\ d\rho}
and
\eq{\label{sec5.3.1f4}V_J^{I,\alpha}(t,x)&:=\int_0^1 (\partial_{v^J_\alpha}Q^I)(\rho \wt{u}+(1-\rho) u,\rho\partial \wt{u}+(1-\rho)\partial u)\ d\rho.}
As a result, $\phi=(\phi^I)$ is a solution to \eq{\label{phieqnsec5}(\wt{\Box}\phi)^I&=\Box \phi^I+\sum_{J=1}^N\sum_{\alpha=0}^3V_J^{I,\alpha} \partial_\alpha \phi^J+\sum_{J=1}^NW_J^I\phi^J=0,\qquad I=1,2,\dots,N.}

Now Theorem part (i) of \ref{mthm0} results from the following lemma. In the proof, we shall make use of both the null condition \eqref{swea2} and the estimate \eqref{def1a2} in Definition \ref{defrad}.

\lem{The functions $V_{J}^{I,\alpha}$, $W^I_J$ and $\phi^J$ defined above satisfy the estimates \eqref{mthma1}-\eqref{mthma4}  in Theorem \ref{mthm}.}
\begin{proof}For each $(u,v)\in\R^{N}\times\R^{4N}$, we have
\fm{\partial_{u^J}Q^I(u,v)=O((|u|+|v|)^2),}
\fm{\partial_{v^J_\alpha}Q^I(u,v)=\sum_{K=1}^N\sum_{\beta=0}^3(A^{\alpha\beta}_{I,JK}+A^{\alpha\beta}_{I,KJ}) v_\beta^K+O((|u|+|v|)^2).}
Note that these estimates from the Taylor expansion \eqref{swea1}.
By \eqref{sec5.3.1f4}, we have
\fm{|V_J^{I,\alpha}|&\leq \int_0^1 |(\partial_{v^J_\alpha}Q^I)(\rho \wt{u}+(1-\rho) u,\rho\partial \wt{u}+(1-\rho)\partial u)|\ d\rho\\
&\lesssim \int_0^1 |\partial u|+|\partial \wt{u}|+(|u|+|\partial u|+|\wt{u}|+|\partial\wt{u}|)^2\ d\rho.}
Note that $\D$ defined in Theorem \ref{mthm} is contained in $\{r-t>R_2>0\}$. Also  recall that we have \fm{|u|+|\partial u|+|\wt{u}|+|\partial\wt{u}|\lesssim (1+t+r)^{-1},\qquad r-R_2>t>0} by \eqref{defradcor} which is a corollary of the definition of radiation fields. As a result, for each $(t,x)\in\D$ we have
\fm{|V_J^{I,\alpha}|&\lesssim (1+t+r)^{-1}+(1+t+r)^{-2}\lesssim \lra{r+t}^{-1}.}
That is, \eqref{mthma1} holds.

Next, by \eqref{sec5.3.1f4}, we have
\eq{\label{lemfff1}&\sum_\alpha V_J^{I,\alpha}\widehat{\omega}_\alpha\\&=\int_0^1 \sum_\alpha\widehat{\omega}_\alpha(\partial_{v^J_\alpha}Q^I)(\rho \wt{u}+(1-\rho) u,\rho\partial \wt{u}+(1-\rho)\partial u)\ d\rho\\
&=\int_0^1 \sum_{K,\alpha,\beta}(A^{\alpha\beta}_{I,JK}+A^{\alpha\beta}_{I,KJ})\widehat{\omega}_\alpha[\rho\partial_\beta \wt{u}^K+(1-\rho)\partial_\beta u^K]\ d\rho+O((|u|+|\partial u|+|\wt{u}|+|\partial\wt{u}|)^2).}
Note that \fm{\partial_\beta u^K=-\widehat{\omega}_\beta (\partial_t-\partial_r)u^K+f_0\cdot(\partial_t+\partial_r)u^K+\sum_{j=1}^3f_0\cdot(\partial_j-\omega_j\partial_r)u^K}
where $f_0$ denotes a polynomial of $\omega$. By \eqref{defradcor} and the finite speed of propagation, we have
\fm{\sum_{\alpha,\beta} A_{I,JK}^{\alpha\beta}\widehat{\omega}_\alpha\partial_\beta u^K&=\sum_{\alpha,\beta} A_{I,JK}^{\alpha\beta}\widehat{\omega}_\alpha\widehat{\omega}_\beta (\partial_t-\partial_r)u^K+O(\lra{t+r}^{-1}|Zu|)=O(\lra{t+r}^{-2}).}
Here we make use of the null condition \eqref{swea2}. We can control the remaining terms in \eqref{lemfff1} by following the same method. As a result, we have $\sum_\alpha V_J^{I,\alpha}\widehat{\omega}_\alpha=O(\lra{t+r}^{-2})$ in $\D$. In addition, by \eqref{sec5.3.1f3}, in $\D$ we have
\fm{|W_J^I|&\lesssim\int_0^1 |(\partial_{u^J}Q^I)(\rho \wt{u}+(1-\rho) u,\rho\partial \wt{u}+(1-\rho)\partial u)|\ d\rho\\
&\lesssim (|u|+|\partial u|+|\wt{u}|+|\partial\wt{u}|)^2\lesssim\lra{t+r}^{-2}.}
As a result, we obtain \eqref{mthma2}.

To prove \eqref{mthma3} and \eqref{mthma4}, we apply Proposition \ref{radfie} with $M=1$. Since $F_0(q,\omega)=\wt{F}_0(q,\omega)$ for all $q> R_2$ and $\omega\in\mathbb{S}^2$, we have
\fm{\sum_{|L|\leq 1}|Z^L\phi|=\sum_{|L|\leq 1}|Z^L(u-\wt{u})|\lesssim  \lra{t}^{-2},\qquad\text{whenever }r-t> R_2.}
The estimate \eqref{mthma3} is obvious, and the estimate \eqref{mthma4} follows from Lemma \ref{l2.1}. 
\end{proof}\rm

\subsubsection{Proof of part (ii) in Theorem \ref{mthm0}} We now assume that $R_2\leq -R_1$. By part (i), we know that $u=\wt{u}$ in $\D_{R_1,R_2'}$ for each $|R_2'|<R_1$. Now we fix $(t,x)\in \R^{1+3}$ such that $t>0$, $r-t>-R_1$ and $r+t>R_1$. Since
\fm{\lim_{R_2'\downarrow -R_1}(r^2-(t+R'_2)^2-(R_1^2-(R_2')^2))=r^2-(t-R_1)^2=(r-t+R_1)(r+t-R_1)>0,}
we can find some $R_2'>R_1$ such that $r^2-(t+R'_2)^2-(R_1^2-(R_2')^2)>0$. In other words, we have
\fm{\{t>0,\ r-t>-R_1,\ r+t>R_1\}\subset \bigcup_{R_2'>-R_1}\D_{R_1,R_2'}}
and as a result, $u=\wt{u}$ whenever $t>0$, $r-t>-R_1$ and $r+t>R_1$.

In particular, we have proved that $(u,u_t)|_{t=R_1}=(\wt{u},\wt{u}_t)|_{t=R_1}$ everywhere except at $x=0$. Since $u$ and $\wt{u}$ are $C^1$ functions, we conclude that $(u,u_t)|_{t=R_1}=(\wt{u},\wt{u}_t)|_{t=R_1}$ everywhere. In other words, the difference $\phi=u-\wt{u}$ is a solution to \eqref{phieqnsec5} with zero data at $t=R_1$. By applying Theorem I.2.2 in Sogge \cite{1sogg}, we conclude that $\phi=0$ everywhere.

\subsubsection{Proof of the result in Remark \ref{rmkk1.4}.}
Let $u$ be a smooth $\R$-valued function such that $\Box u=0$. Suppose that $u=0$ whenever $t>0$ and $r^2-(t+R_2)^2>R_1^2-R_2^2$ for some $R_1>0$ and $|R_2|<R_1$. We  claim that $u=0$ whenever $r+t>R_1$ and $r-t>R_2$.

Our main tool is the Holmgren's theorem. See the discussion in Section \ref{sec1.1.1}. We define a new function
\fm{\psi_\kappa(t,x):=(r-\frac{R_1+R_2}{2}-\kappa)^2-(t-\frac{R_1-R_2}{2})^2.}
Here $\kappa>0$ is a small constant.

We first check that the level set $\{\psi_\kappa=c\}$ for a fixed constant $c>0$ is noncharacteristic everywhere with respect to $\Box$. The principal symbol of $\Box$ is $p(\xi)=-m^{\alpha\beta}\xi_\alpha\xi_\beta$, and we have \fm{\partial_t\psi_\kappa=-2(t-\frac{R_1-R_2}{2}),\qquad \partial_j\psi_\kappa=2(r-\frac{R_1+R_2}{2}-\kappa)\omega_j.} Then, we have
\fm{p(\nabla_{t,x}\psi_\kappa)&=[-2(t-\frac{R_1-R_2}{2})]^2-\sum_{j=1}^3[2(r-\frac{R_1+R_2}{2}-\kappa)\omega_j]^2=-4\psi_\kappa=-4c<0.}
As a result, the level set $\{\psi=c\}$ is noncharacteristic whenever $c>0$.

Moreover, for each fixed $c_0>0$, the set
\fm{\bigcup_{c>c_0}\{(t,x):\ t>0,\ \psi_\kappa(t,x)=c\}\setminus\D_{R_1,R_2}}
is a bounded set in $[T_1,T_2]\times\R^3$ for some $0<T_1<T_2<\infty$. Then, we can follow the proof in Section \ref{sec5.2} to prove that $u=0$ in
\fm{\bigcup_{c>c_0}\{(t,x):\ t>0,\ \psi_\kappa(t,x)=c\}}for each fixed $c_0>0$ and $\kappa>0$. It is easy to show that
\fm{\bigcup_{c_0>0}\bigcup_{\kappa>0}\bigcup_{c>c_0}\{(t,x):\ t>0,\ \psi_\kappa(t,x)=c\}=\{(t,x):\ t>0,\ r-t>R_2,\ r+t>R_1\}.}
This finishes the proof.\rm

\bibliography{paperucp}{}
\bibliographystyle{plain}

\end{document}